\documentclass[10pt]{article}

\usepackage{bbold}
\usepackage{color}
\usepackage{graphics} % for pdf, bitmapped graphics files
\usepackage{needspace}

 \input{mysymbol.sty}
\usepackage{cite}
\usepackage{amsmath}
\usepackage{amssymb}
\usepackage{amsthm}
\usepackage{mathtools}
\usepackage{multirow}
\usepackage{url}
\usepackage{graphics, subfigure, times, amsfonts}
\usepackage{tikz, epic,eepic}
\usetikzlibrary{shapes,arrows}
\usepackage{pgfplots}
\usepackage{color}
\usepackage{hyperref}
\usepackage{epstopdf}
\usepackage{epsfig, amsbsy}
\usepackage{latexsym}
\usepackage{amscd, verbatim}
\usepackage{multirow}
\usepackage{booktabs}

\newtheorem{theorem}{Theorem}
\newtheorem{claim}{Claim}
\newtheorem{proposition}{Proposition}
\newtheorem{corollary}{Corollary}
\newtheorem{lemma}{Lemma}
\newtheorem{definition}{Definition}

{\itshape}{\rmfamily}

\newtheorem{remark}{Remark}

\usepackage{tikz}
\usetikzlibrary{shapes.geometric, arrows}
\tikzstyle{area} = [rectangle, rounded corners, minimum width=8cm, minimum height=3cm,text centered, text width=7cm,draw=black, fill=white!40]
\tikzstyle{arrow} = [thick,->,>=stealth]

\usepackage[ruled,vlined]{algorithm2e}   

\usepackage{textcomp}

\title{\Large An Iterative Mechanism for Coupling Electricity Markets }

\author{Alfredo Garcia, Roohallah Khatami, Ceyhun Eksin and Furkan Sezer
\thanks{
A. Garcia,  C. Eksin and F. Sezer are with the Department of Industrial and Systems Engineering, Texas A\&M University, College Station, TX 77843 USA (e-mails: alfredo.garcia@tamu.edu, eksinc@tamu.edu, furkan.sezer@tamu.edu)
}}

\def\bm#1{\mbox{\boldmath $#1$}}

\newcommand{\beq}{\begin{equation}}
\newcommand{\eeq}{\end{equation}}
\newcommand{\bea}{\begin{array}}
\newcommand{\ena}{\end{array}}
\newcommand{\bds}{\begin {itemize}}
\newcommand{\eds}{\end {itemize}}
\newcommand{\bdf}{\begin{definition}}
\newcommand{\edf}{\end{definition}}
\newcommand{\blm}{\begin{lemma}}
\newcommand{\elm}{\end{lemma}}
\newcommand{\bthm}{\begin{theorem}}
\newcommand{\ethm}{\end{theorem}}
\newcommand{\bprp}{\begin{prop}}
\newcommand{\eprp}{\end{prop}}
\newcommand{\bcl}{\begin{claim}}
\newcommand{\ecl}{\end{claim}}
\newcommand{\bcr}{\begin{coro}}
\newcommand{\ecr}{\end{coro}}
\newcommand{\bquest}{\begin{question}}
\newcommand{\equest}{\end{question}}

\newcommand{\larrow}{{\larrow}}

\begin{document}

\maketitle

\begin{abstract}
The coordinated operation of interconnected but locally controlled electricity markets is generally referred to as a ``coupling". In this paper we propose a new mechanism design for efficient coupling of independent electricity markets.
The mechanism operates after each individual market has settled (e.g. hour-ahead) and based upon the reported supply and demand functions for internal market optimization (clearing), each market operator is asked to iteratively quote the terms of energy trade (on behalf of the agents participating in its market) across the transmission lines connecting to other markets.
The mechanism is scalable as the informational demands placed on each market operator at each iteration are limited. We show that the mechanism's outcome converges to the optimal flows between markets given the reported supply and demand functions from each individual market clearing.  We show the proposed market coupling design does not alter the structure of incentives in each internal market, i.e., any internal market equilibrium will remain so (approximately) after coupling is implemented. This is achieved via incentive transfers (updated at each iteration) that remunerate each market with its marginal contribution (i.e. cost savings) to all other participating markets.
We identify a sufficient condition on a uniform participation fee for each market operator ensuring the mechanism incurs no deficit.
The proposed decentralized mechanism is implemented on the three-area IEEE Reliability Test System where the simulation results showcase the efficiency of proposed model.

\end{abstract}

%\begin{IEEEkeywords}
%distributed optimization, subgradient algorithms
%\end{IEEEkeywords}

% \IEEEpeerreviewmaketitle

\section{Introduction}

The control and operation of many interconnected electricity markets (e.g. Western Europe, Northeast US) is conducted in a decentralized manner by independent agencies (e.g. system operators ISOs in the US, power exchanges PXs and transmission system operators TSOs in Europe). In a decentralized market architecture, each individual market operator has the necessary control authority over an area to ensure the associated power system operation achieves an acceptable trade-off between economic efficiency and reliability \cite{Meyn, molzahn2017survey}. However, the interconnection of decentralized electricity markets by means of transmission lines (or tielines) provides increased efficiency and/or reliability by enabling access to cheaper generation and/or flexibility in the form of reserves. The coordinated operation of interconnected but locally controlled electricity markets is generally referred to as market coupling. Such coupling can be seen as a consensus amongst independent market operators on the technical and economic terms associated to energy flows between markets.

In this paper we present a market coupling mechanism that operates after each individual market has settled (e.g. day-ahead hourly).  Based upon the reported supply and demand cost functions, each market operator is asked to iteratively quote the terms of energy trade (quantities and prices) across interties (on behalf of the agents participating in its market) by recomputing internal market clearing dispatch, a task that involves solving a DC optimal power flow (a convex optimization problem). Each market operator exchanges voltage phase angles of boundary buses with adjacent market operators. This information is needed to ensure the flows across markets are feasible.
The market coupling design is {\em iterative}: a capacity price for each intertie is updated at every round (or iteration) to reflect excess demand and  each market operator exchanges voltage phase angles of boundary buses with adjacent area operators (Section \ref{sec.alloc}). We show that the sequence of outcomes of the proposed market coupling mechanism converges to the optimal allocation and pricing of intertie capacity {\em given} the reported supply and demand functions from each individual market clearing (Theorem \ref{thm_convergence}). In other words, the limit outcome corresponds to the solution of the joint (centralized) DC-OPF for all areas.

In the second part of the paper we focus on incentive compatibility issues that may {\em exclusively} arise from the pricing and allocation rules of intertie capacity (Section \ref{sec.inc}). To focus on incentive compatibility issues related to market coupling, we assume the reported supply and demand functions for the clearing of each market (e.g. hour-ahead) are a {\em Nash equilibrium} and examine the ways in which such equilibrium may be distorted. Since the net surplus gains resulting from coupling are to be {\em fully distributed} back to market participants, any agent that attempts a manipulation (by deviating from equilibrium reporting) receives a {\em fraction} of surplus gains. Therefore, in order to show that there is no incentive to manipulate the reported supply or demand functions it is sufficient to demonstrate that any such manipulation does not yield increased net surplus at the {\em aggregate} market level. Indeed, we show that if supply and demand cost functions for internal market clearing are a Nash equilibrium, they remain an approximate Nash equilibrium after the coupling of areas is implemented (Theorem \ref{thm_Nash}). Thus the proposed market coupling design has {\em negligible effect} on the structure of competition in each individual market. This result relies on incentive (monetary) transfers to each market operator which are set to equal an estimate of each market's {\em marginal} contribution to the coupling, i.e. 
the change in surplus that a given area's participation induces in all other areas.
Such estimate is a function of reported information by market operators (i.e. desired intertie flows and locational marginal prices at boundary nodes). Setting incentive transfers to approximately equal each individual's market marginal contribution to the coupling, aligns the incentives of each individual market with that of maximizing surplus gains from trade\footnote{ Or equivalently, minimizing total cost when demand is inelastic.} net of any congestion rent associated with intertie trades.

To ensure the mechanism does not incur a deficit, each area is assessed a participation fee which must be bounded below by the value of the lowest individual surplus gain from coupling across areas.
We identify a sufficient condition for ensuring no deficit which corresponds to the case in which the {\em information rents}, i.e. the costs incurred to ensure incentive compatibility, do not exceed the system-wide benefits obtained from coupling (Theorem \ref{thm_balance}).
Numerical experiments on a three-zone IEEE Reliability Test System are used to illustrate our results (Section \ref{sec.results}). 

% {\em only reported information} is used during the implementation of the mechanism.

% that computes incentive transfers and participation fees based on locally available information 

% In section \ref{practical}, we show these results continue to hold when {\em only reported information} is used during the implementation of the mechanism.
%  In section \ref{sec.results}, we provide a numerical illustration with a three-zone IEEE Reliability Test System. We end the paper with conclusions.
 
 \section{Literature Review}
 
The market coupling problem is related to the literature in distributed consensus optimization algorithms for resource allocation with coupled constraints (see for example \cite{Chang} and \cite{Prandini}). However, this literature does not address incentive compatibility concerns arising from potential strategic behavior by participating agents, an important consideration for coupling electricity markets. The iterative mechanism is akin to best-response type Nash equilibrium seeking algorithms in which updates are reported with inertia \cite{Shamma_Arslan_2005,marden2009joint,swenson2018distributed,grammatico2017dynamic,parise2019distributed}. Here again in addition to convergence we provide remuneration schemes that make sure internal equilibrium conditions in each market are not altered as a result of coupling.

The paper is also related to the literature on mechanism design on networks \cite{Barrera2015}, \cite{Farhadi} and \cite{Sinha} which studies incentives for efficient allocation of resources to agents interacting in a network. In contrast, in this paper we focus on incentives for the efficient interconnection of multiple networks (i.e. markets) via interties. Our overall goal is different as we seek to identify pricing and allocation rules for interconnection capacity that are guaranteed to not distort internal equilibrium conditions in each market and enable the identification of optimal allocation and pricing of intertie capacity {\em given} equilibrium supply and demand functions from each individual market clearing.

There is also a related literature on distributed solutions to the economic dispatch problem in power systems to increase scalability, robustness to failure and cyber-security \cite{molzahn2017survey, kargarian2016toward,Cherukuri_2016, Cherukuri_2020, kar2012distributed,dominguez2012decentralized,yang2013consensus,binetti2014distributed,dorfler2015breaking,cherukuri2015distributed}. In contrast, in this paper we focus on coupling the operation of multiple independently operated markets via interties wherein each market optimizes its own internal operation but consensus must be achieved on the flows across markets. 

Finally, the iterative mechanism considered in this paper is related to decentralized approaches for solving the joint interconnected economic dispatch problem in a distributed fashion via primal decomposition methods \cite{bakirtzis2003decentralized}, \cite{Zhao}, \cite{Li}, \cite{Guo_2} or dual decomposition methods \cite{Conejo}, \cite{binetti2014distributed}, \cite{erseghe2014distributed}. Improvements to existing market coupling designs are considered in \cite{Guo}, \cite{Ji}. However, this literature does not address incentive compatibility concerns that may be induced by market coupling.

The European Price Coupling of Regions (PCR) project is perhaps the largest example of a market coupling mechanism currently in operation. In such mechanism, on a day-ahead basis market participants submit their orders to their respective power exchange (PX). These orders are collected and submitted to a centralized market-clearing algorithm ({\em Euphemia}) which identifies the prices and the trades that should be executed so as to maximize the gains from trade generated by the executed orders while ensuring the available transfer capacity is not exceeded \cite{ENTSO-E, lam2018european}.

In the United States, a different approach referred to as coordinated transaction scheduling (CTS) has been approved by FERC and implemented in NYISO,\footnote{https://www.nyiso.com/documents/20142/3036853/CTS+ISONE.pdf} PJM,\footnote{https://www.pjm.com/-/media/committees-groups/stakeholder-meetings/pjm-nyiso/20180402/20180402-item-02b-coordinated-transaction-scheduling-metrics.ashx} and MISO.\footnote{https://www.misoenergy.org/stakeholder-engagement/issue-tracking/miso-pjm-coordinated-transaction-scheduling-cts/}
In this method, proxy buses at the interties are used as trading points for enabling bids from external market participants. Thus offers to sell (buy) from external market participants are represented by injections (withdrawals) at such proxy buses.

There are important concerns regarding both the PCR and CTS designs.
First, participants may exercise market power by reporting information (e.g. orders to buy or sell) that is not consistent with the overall goal of maximizing social surplus. That is, the schemes are not guaranteed to be {\em incentive compatible}. Thus, for example, with distorted information as input, it does not follow that a sophisticated market-clearing algorithm such as {\em Euphemia} results in the efficient allocation and pricing of available interconnection capacity.
A second concern pertains to market architecture. The PCR design implies individual area operators play a largely passive role as conveyors of bidding information. This goes against the widely accepted notion that with increasing intermittency (due to increased renewable capacity) individual area operators are likely to have superior information regarding acceptable trade-offs in cost vs security for the operation of their respective areas. A uniform protocol such as PCR limits the discretion of area operators in implementing heterogeneous approaches for market design and operation. 
CTS is limited to scheduling interchanges between two neighboring areas. This clearly limits efficiency gains as areas may need to interact with more than one neighboring area simultaneously in order to identify maximum gains from trade. In addition, the fictitious nature of proxy buses may imply that the actual power flow may differ from the CTS schedule. Thus, it is of interest to identify alternative market coupling designs ensuring efficiency with incentive compatibility guarantees while giving a more active role to individual area operators.

\section{Tieline Capacity Allocation and Pricing} \label{sec.iter}
\begin{figure}[h] 
	\centering
\vspace{-5pt}
	\includegraphics[width=0.5\linewidth]{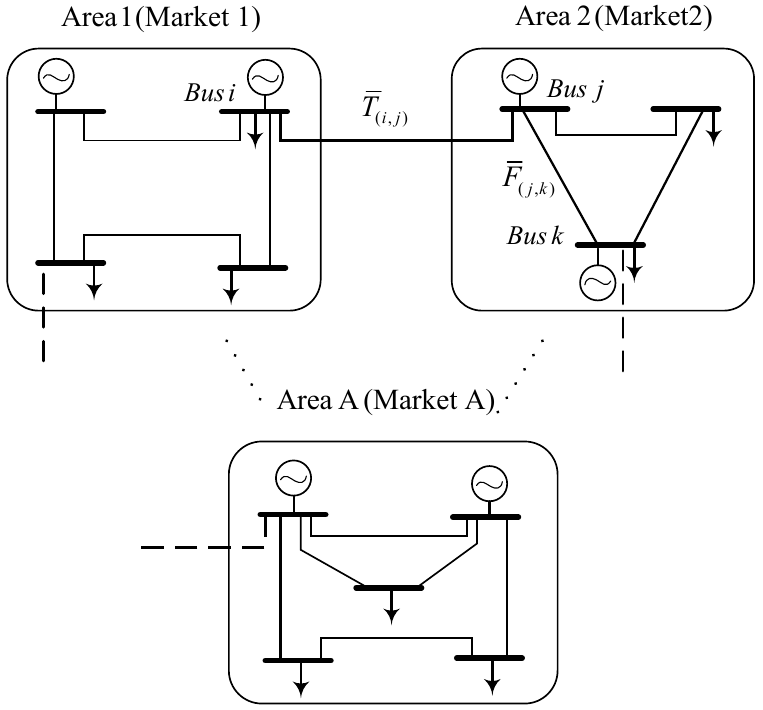}
     \vspace{-5pt}
	\caption{\small Generic multi-area power transmission network.}
       \label{fig.network}
\end{figure}

Consider a power system composed of $A>1$ operational areas, as schematically shown in Fig. \ref{fig.network}, where the areas are corresponded to the associated electricity markets governed by separate ISOs\footnote{In the rest of this paper we use the terms ``areas", ``markets", and ``ISOs" interchangeably.}. 
Each area $a\in\mathcal A=\{1,2,\ldots,A\}$ is represented by a directed graph $(\mathcal{N}_a,\mathcal{F}_a)$ where $\mathcal{N}_a=\{1,2,\ldots,N_a\}$ and $\mathcal{F}_a=\{(i,j)|i,j \in \mathcal{N}_a,j\equiv j(i)\}$ respectively represent the set of nodes (buses) and edges (transmission lines). 
For area $a \in \mathcal{A}$, the bus voltage angles and nodal loads form respectively the vectors $\bm{\theta}_a=\left(\theta_{a,1},\theta_{a,2},\ldots,\theta_{a,N_a}\right)^T$ and $\mathbf{D}_a=\left(D_{a,1},D_{a,2},\ldots,D_{a,N_a}\right)^T$. The $N_a \times N_a$ admittance matrix of each area is denoted by $\mathbf{B}_a$.
The set of $G_a$ generating units at each area is represented by $\mathcal{G}_a=\{1,\ldots,G_a\}$, the power generation of units form the vector $\mathbf{P}_a=\left(P_{a,1},P_{a,2},\ldots,P_{a,G_a}\right)^T$, and the $N_a \times G_a$ incidence matrix $\mathbf{M}_a$ maps the generating units to buses. 
The tielines entering/leaving  area $a$ are defined by set $\mathcal{T}_a=\{(i,j)|i\in \mathcal{N}_a, j\in \mathcal{N}_{a'}, a'\in \mathcal A\backslash a \}$, the associated tieline power flows form the vector $\mathbf{T}_a=\left(T_{a,(i,j)}\right),~(i,j)\in \mathcal{T}_a$ , and the incidence matrix mapping tieline power flows to buses is shown by $\mathbf{R}_a$. 

\subsection{Centralized DC-OPF} \label{sec.CentOPF}
In describing the operation of the proposed market coupling mechanism, it will be useful to refer to the centralized DC-OPF problem where  a single system operator optimally operates the entire power system with complete information regarding all the areas' generators supply costs. \footnote{For ease of exposition, we leave out demand offers and assume demand is inelastic and known. This implies the internal market clearing corresponds to solving the DC power flow that minimizes total cost. The results can be extended with active demand participation and market clearing modeled as DC power flow maximizing aggregate surplus.} This is a hypothetical exercise as in practice there is no single market operator and information on generators costs is not available to any single party.
In the formulation below, the areas are separated and the transmission line power flows are distinguished from tieline power flows \cite{bakirtzis2003decentralized}.
This arrangement would facilitate formulating the decentralized DC-OPF, which is the backbone of proposed tieline capacity allocation.

To sum up, in the centralized DC-OPF formulation the total operational cost of power system in \eqref{eq.Obj1} is minimized subject to the operation constraints:
\begin{align} \label{eq.Obj1}
\min_{P_{a,g}} &~~~\sum_{a\in \mathcal A} \sum_{g\in \mathcal{G}_a} C_{a,g}(P_{a,g}),\\  
\label{eq.Bal1}
&\mathbf{B}_a\bm{\theta}_a + \mathbf{R}_a \mathbf{T}_a=\mathbf{M}_a \mathbf{P}_a - \mathbf{D}_a, \quad ({\boldsymbol{\alpha }}_{a}), ~~\forall a \in \mathcal A,\\ \label{eq.Plim1}
&\underline{\mathbf{P}}_a\leq \mathbf{P}_a\leq  \overline{\mathbf{P}}_a, \quad  (\nu _{a,g},\lambda
_{a,g}), ~~ \forall a \in \mathcal A,\\ \label{eq.Flim1}
-&\overline{F}_{a,(i,j)}\!\leq\! \frac{\theta_{a,i}-\theta_{a,j}}{x_{a,(i,j)}} \!\leq\! \overline{F}_{a,(i,j)},\nonumber \\& \qquad \qquad (\kappa _{a,(i,j)},\eta
_{a,(i,j)}),  ~~\forall a \in \mathcal A, \forall (i,j) \in \mathcal{F}_a,\\  \label{eq.TieFlow1}
&T_{a,(i,j)}= \frac{\theta_{a,i}-\theta_{a',j}}{\bar{x}_{a,(i,j)}}, \nonumber\\&\qquad \qquad  \quad(\xi _{a,(i,j)}),~~\forall a \in \mathcal A,~\forall  (i,j) \in \mathcal{T}_a,\\  \label{eq.Tielim1}
-&\overline{T}_{a,(i,j)}\!\leq\! T_{a,(i,j)}\!\leq\! \overline{T}_{a,(i,j)}, \quad  \nonumber \\& \qquad \qquad (\bar{\kappa}_{a,(i,j)},\bar{\eta}_{a,(i,j)}), ~\forall a \in \mathcal A,~\forall  (i,j) \in \mathcal{T}_a,\\  \label{eq.slack1}
&\theta_{1,1}=0,
\end{align} 
where $C_{a,g}(P_{a,g})$ denotes the strictly convex cost functions of generating units, the vectors $\underline{\mathbf{P}}_a$ and $\overline{\mathbf{P}}_a$ respectively represent the minimum and maximum generation limits of units, and $\overline{F}_{a,(i,j)}$ and $\overline{T}_{a,(i,j)}$ respectively represent the power flow limits of transmission lines and tielines.
The nodal power balance is secured through \eqref{eq.Bal1}, the generation of units are confined to their limits in \eqref{eq.Plim1}, transmission line power flows are maintained within thermal limits in \eqref{eq.Flim1} where $x_{a,(i,j)}$ represents the reactance of each line, and the tieline power flows are calculated and constrained to their thermal limits in \eqref{eq.TieFlow1} and \eqref{eq.Tielim1} where $\bar{x}_{a,(i,j)}$ represents the tieline reactances.
Further, the voltage phase angle of the slack bus, which is assumed to be the bus 1 of area 1, is set to zero in \eqref{eq.slack1}.
Note that $\overline{T}_{a,(i,j)}=\overline{T}_{a',(j,i)}$ and $\bar{x}_{a,(i,j)}=\bar{x}_{a',(j,i)}$. 

The first order conditions with respect to $\theta _{a,i}$, $\theta _{a^{\prime },j}$%
and $T_{a,(i,j)}$ for tieline $(i,j)\in \mathcal{T}_{a}$  in the centralized DC-OPF problem   \eqref {eq.Obj1}- \eqref{eq.slack1} are:%
\begin{align}\label{eq.centralizedangle_1}
\sum\limits_{\left. j\right\vert (i,j)\in \mathcal{F}_{a}}\Big[\frac{1}{x_{a,(i,j)}}(\alpha^{\ast}
_{a,i}-\alpha^{\ast}_{a,j}&+\eta^{\ast}_{a,(i,j)} -\kappa^{\ast} _{a,(i,j)}) -\frac{1}{\bar{x}_{a,(i,j)}}\xi^{\ast}
_{a,(i,j)}\Big]=0
\end{align}

\begin{align} \label{eq.centralizedangle_2}
\sum\limits_{\left. j\right\vert (i,j)\in \mathcal{F}_{a^{\prime
}}}\Big[\frac{1}{x_{a',(i,j)}}(\alpha^{\ast} _{a^{\prime },i}-&\alpha^{\ast} _{a^{\prime },j}+\eta^{\ast} _{a^{\prime
},(i,j)}-\kappa^{\ast} _{a^{\prime },(i,j)})+\frac{1}{\bar{x}_{a,(i,j)}}\xi^{\ast} _{a,(i,j)}\Big]=0
\end{align}

\begin{equation} \label{eq_centralized_intertie}
\alpha^{\ast} _{a,i}-\alpha^{\ast} _{a^{\prime },j}+\bar{\eta}^{\ast} _{a,(i,j)}-\bar{\kappa}^{\ast}
_{a,(i,j)}+\xi^{\ast} _{a,(i,j)}=0
\end{equation}
The centralized model presented in this section leads to an efficient tieline capacity allocation, yet it requires the areas to reveal their private information (network configuration and cost functions of units) and leads to solving an optimization problem of high dimensionality.  
However, the decentralized DC-OPF formulation and the counterpart iterative capacity allocation method presented in the following sections reduce the size of optimization problems and preserve the information privacy of areas.

\subsection{Decentralized DC-OPF} \label{sec.DisOPF}
In this section, we introduce a decentralized approach for solving the DC-OPF problem. Given locational marginal prices ($\alpha_{a',j}$) and angles ($\theta_{a',j}$) of adjacent areas and the capacity price $\mu_{(i,j)}$, the problem for area $a \in \mathcal{A}$ is given by:

% \begin{align} \label{eq.Obj2}
% &\min ~~ V_a = \sum_{g\in \mathcal{G}_a} C_{a,g}(P_{a,g}) \nonumber\\
% &~~~~~~~~~~~+\!\!\!\!\!\sum_{(i,j)\in \mathcal{T}_a} \!\!\!\big[\frac{\mu_{i,j}}{2}|T_{a,(i,j)}|\!-\!\alpha_{a',(i,j)} T_{a,{(i,j)}}\big],\\ 
% \label{eq.Bal2}
% &\mathbf{B}_a\bm{\theta}_a + \mathbf{R}_a \mathbf{T}_a=\mathbf{M}_a \mathbf{P}_a - \mathbf{D}_a,~~~~\left({\boldsymbol\alpha}_a\right),\\
% \label{eq.Plim2}
% &\underline{\mathbf{P}}_a\leq \mathbf{P}_a\leq  \overline{\mathbf{P}}_a,~~~~~~~~~~~~~~~~~~\left(\underline{\boldsymbol{\nu}}_a,\overline{\boldsymbol{\nu}}_a\right),\\
% \label{eq.Flim2}
% -&\overline{F}_{a,(i,j)}\!\leq\! \frac{\theta_{a,i}-\theta_{a,j}}{x_{a,(i,j)}} \!\leq\! \overline{F}_{a,(i,j)}, \forall (i,j) \in \mathcal{F}_a,\nonumber\\
% &\left(\underline{\eta}_{a,(i,j)},\overline{\eta}_{a,(i,j)}\right)\\
%  \label{eq.TieFlow2}
% &T^k_{a,(i,j)}= \frac{\theta_{a,i}-\theta_{a',j}}{x^t_{a,(i,j)}}, ~~\forall  (i,j) \in \mathcal{T}_a,\left(\underline{\zeta}_{a,(i,j)}\right),\\
% \label{eq.slack2}
% &\theta_{a,1}=0, ~ \mbox{if} ~ a=1,
% \end{align}

\begin{align} \label{eq.Obj2}
\min &~~~ \sum_{g\in \mathcal{G}_a} C_{a,g}(P_{a,g}) -\sum_{(i,j)\in \mathcal{T}_a} \alpha_{a',j} T_{a,{(i,j)}}+\sum_{(i,j)\in \mathcal{T}_a}\frac{\mu_{(i,j)}}{2}(|T_{a,(i,j)}|-\bar{T}_{a,(i,j)})
\\ 
\label{eq.Bal2}
&\mathbf{B}_a\bm{\theta}_a + \mathbf{R}_a \mathbf{T}_a=\mathbf{M}_a \mathbf{P}_a - \mathbf{D}_a,~~~~\left(\hat{\boldsymbol\alpha}_a\right),\\
\label{eq.Plim2}
&\underline{\mathbf{P}}_a\leq \mathbf{P}_a\leq  \overline{\mathbf{P}}_a,\\
\label{eq.Flim2}
-&\overline{F}_{a,(i,j)}\!\leq\! \frac{\theta_{a,i}-\theta_{a,j}}{x_{a,(i,j)}} \!\leq\! \overline{F}_{a,(i,j)},\forall (i,j) \in \mathcal{F}_a, \quad(\hat{\kappa} _{a,(i,j)},\hat{\eta}_{a,(i,j)}) \\
 \label{eq.TieFlow2}
&T_{a,(i,j)}= \frac{\theta_{a,i}-\theta_{a',j}}{\bar{x}_{a,(i,j)}}, ~~\forall  (i,j) \in \mathcal{T}_a,\left(\hat{\xi}_{a,(i,j)}\right),\\
\label{eq.slack2}
&\theta_{a,1}=0, ~ \mbox{if} ~ a=1.
\end{align}
Note the objective function in \eqref{eq.Obj2} not only includes area $a$'s total generation cost but also accounts for the capacity costs of incident tielines as well as the cost/revenue of the energy import/export from/to neighbouring areas.
It is assumed that the interconnected areas split the capacity cost over tielines, meaning that the two incident areas of a tieline are both charged for the excess capacity usage at the rate of $\frac{\mu_{(i,j)}}{2}$. This is done with no loss of generality: other conventions (e.g. exporter is assigned capacity cost, importer is assigned capacity cost) can also be used without altering the results in this paper.

The following lemma makes precise in which way decentralization is made possible in this formulation.

\begin{lemma} \label{lemma_equivalence}
If the locational marginal prices, angles of adjacent areas and the capacity price of interties correspond to the optimal solution of centralized OPF (i.e. $\alpha_{a',j}^{*}$, $\theta_{a',j}^{*}$ and $\frac{\mu^{\ast}_{(i,j)}}{2}\mbox{sign} ({T}_{a,(i,j)}^*)=\bar{\eta}^{\ast}_{a,(i,j)}-\bar{\kappa}^{\ast}_{a,(i,j)}$), then the solution to the decentralized OPF problem for area $a \in \mathcal{A}$ corresponds to the solution of the centralized  OPF.

\end{lemma}

\begin{proof}
Consider  the first order conditions with respect to $\theta _{a,i}$ and $T_{a,(i,j)}$ in the decentralized DC-OPF problem \eqref{eq.Obj2}-\eqref{eq.slack2} for area $a \in \mathcal{A}$:
\begin{align}
\sum\limits_{\left. j\right\vert (i,j)\in \mathcal{F}_{a}}\Big[\frac{1}{x_{a,(i,j)}}(\hat{\alpha}
_{a,i}-\hat{\alpha}_{a,j}&+\hat{\eta}_{a,(i,j)}-\hat{\kappa}_{a,(i,j)})-\frac{1}{\bar{x}_{a,(i,j)}}\hat{\xi}
_{a,(i,j)}\Big]=0
\end{align}%
\begin{equation}
\hat{\alpha}_{a,i}+\hat{\xi}
_{a,(i,j)}=\alpha^{\ast} _{a^{\prime },j} -\frac{\mu^{\ast}_{(i,j)}}{2} \mbox{sign} (\hat{T}_{a,(i,j)})
\end{equation}
where $\hat{T}_{a,(i,j)}$ is the optimal transfer flow for area $a$ and $\mbox{sign}(x)=1$ if $x\geq0$ and $\mbox{sign}(x)=-1$ if $x<0$.

We check now that the centralized DC-OPF solution is a solution to each decentralized DC-OPF problem, namely that conditions (17)-(18) are satisfied with $\hat{\alpha}_{a,i}=\alpha^{\ast}_{a,i}$,
$\hat{\xi}_{a,(i,j)}=\xi^{\ast}_{a,(i,j)}$ and
\begin{equation*}
\hat{T}_{a,(i,j)}=\left\{ 
\begin{array}{rr}
T_{a,(i,j)}^{\ast } & T_{a,(i,j)}^{\ast }\geq 0 \\ 
&  \\ 
-T_{a,(i,j)}^{\ast } & T_{a,(i,j)}^{\ast }<0%
\end{array}%
\right. 
\end{equation*}
Thus, it can be seen that equation (8) is equivalent to (17). Similarly, equation (18) corresponds to equation (10) above.  The remaining optimality conditions related to internal operation of area $a$ match those of the centralized OPF problem.
\end{proof}

\begin{remark} \label{Rem.Consensus}
The decentralization approach in \cite{bakirtzis2003decentralized} has each area internalizing the tieline capacity constraint. However, adjacent areas have no means of ensuring they agree on the marginal value of intertie capacity.
Hence, one can not conclude that the first order conditions for the decentralized problems match those of the centralized DC-OPF. The decentralization approach in \cite{bakirtzis2003decentralized} only works if 
interties have no capacity contraints (so that marginal value of capacity $\mu_{(i,j)}$ is null)
or if an additional consensus constraint of the form
\begin{equation*}
\bar{\eta}_{a,(i,j)}-\bar{\kappa}_{a,(i,j)} + \bar{\eta}
_{a^{\prime },(i,j)}-\bar{\kappa}_{a^{\prime },(i,j)}=0
\end{equation*}
is introduced where $\bar{\eta}_{a,(i,j)}$ (respectively, $\bar{\eta}_{a^{\prime },(i,j)}$) is the Lagrange multiplier associated with upper bounds on intertie flow from $a$ to $a'$ (respectively, from $a'$ to $a$) and $\bar{\kappa}_{a,(i,j)}$ (respectively, $\bar{\kappa}_{a^{\prime},(i,j)}$) is the multiplier associated with lower bound on intertie flow.

\end{remark}

%In addition, as obvious from data flow diagram of Fig. \ref{fig.DataFlow}, the terms with superscript $(k-1)$ are calculated from the previous iteration and considered as parameters in the present iteration, while the terms with superscript $k$ are deemed as decision variables of the present iteration.

%\noindent\textbf{Note:} The import/export bids are the locational marginal prices (LMPs) at the boundary buses.
%In other words, the optimal Lagrange multipliers of the power balance constraint \eqref{eq.Bal2} at the boundary buses represent the import/export bids of the areas for the energy transacted on tielines.

\subsection{Iterative Mechanism for Tieline Capacity Allocation and Pricing} \label{sec.alloc}
In this section we elaborate the proposed capacity allocation and pricing model.
The market coupling design operates {\em after} each individual market has settled (e.g. hour-ahead), thus providing the initial conditions for the mechanism. Then, based upon the reported supply and demand functions each market operator is asked to iteratively quote the terms of energy trade across interties on behalf of the agents participating in its market.

The basic steps of the mechanism are as follows---see Algorithm \ref{al.Iterative} for details. At each iteration $k>0$ of the mechanism,  
\begin{enumerate}
    \item (Information exchange) Each area reports the terms of trade for each interconnection in the form of intertie flows ($T_{a',(i,j)}^{k-1}$), locational marginal prices ($\alpha^{k-1}_{a',j}$) and angles of adjacent areas ($\theta^{k-1}_{a',j}$), where $\alpha^{0}_{a',j}=0$.
\item (Updating Intertie Flows) Given a capacity price $\mu^{k-1}_{i,j}$, each area $a \in \mathcal{A}$ solves \eqref{eq.Obj2} -  \eqref{eq.slack2} and reports the desired flows $\mathbf{\hat{T}}_{a}^{k}$, corresponding angles at interties  $\boldsymbol{\hat{\theta}}_{a}^{k}$ and locational marginal prices $\boldsymbol{\hat{\alpha}}_{a}^{k}$ to the agent in charge of executing the mechanism which updates the flows, angles and locational marginal prices along the interties according to:
\begin{align}
\mathbf{T}_{a}^{k} = & (1-\rho _{k})\mathbf{T}_{a}^{k-1}+\rho _{k} \hat{\mathbf{T}}_{a}^{k} \label{eq_inertia_a}\\
\boldsymbol{\theta}_{a}^{k} = &(1-\rho _{k})\boldsymbol{\theta}_{a}^{k-1}+\rho _{k} 
\boldsymbol{\hat{\theta}}_{a}^{k} \label{eq_inertia_b}\\
\boldsymbol{\alpha}_{a}^{k} = &(1-\rho _{k})\boldsymbol{\alpha}_{a}^{k-1}+\rho _{k} 
\boldsymbol{\hat{\alpha}}_{a}^{k}\label{eq_inertia_c}
\end{align}
with $\rho_k \rightarrow 0^{+}$, $\sum 
_k \rho_k = + \infty $ and $\sum 
_k \rho^{2}_k < \infty$.
\item (Updating Intertie Capacity Prices) 
\begin{align} \label{eq.PriceUpdate}
&\mu _{(i,j)}^{k}=\nonumber \\
&\max \{\mu _{(i,j)}^{k-1}+\beta (\frac{\left\vert
T_{a,(i,j)}^{k}\right\vert +\left\vert T_{a^{\prime },(i,j)}^{k}\right\vert 
}{2}-\bar{T}_{a,(i,j)}),0\}
\end{align}
where $\beta \in (0,1)$ and $\mu _{(i,j)}^{0}= 0$. 

%\red{\{RK: In the simulation results the $\mu^0_{(i,j)}$ is equal to 50. In the initial versions of this work we were assuming $\mu^0_{(i,j)}$ to be a very big number. That is why I chose 50 which is greater than the highest incremental cost rate of the most expensive unit.\}}. 

\item A monetary transfer $\widehat{\Delta \pi}_{a,k}$ is allocated to each area $a \in \mathcal{A}$. Upon stopping after $T>0$ iterations a coupling participation fee $R>0$ is assessed to each area $a \in \mathcal{A}$.

\end{enumerate}

% Figure \ref{fig.DataFlow} summarizes the key update steps within an area and the information exchange between two coupled areas. 

In updating the intertie flows and angles, and locational marginal prices, each area uses an inertial update in \eqref{eq_inertia_a}-\eqref{eq_inertia_c} by weighting the current optimal solution, i.e., $\mathbf{\hat{T}}_{a}^{k}$, $\boldsymbol{\hat{\theta}}_{a}^{k}$, $\boldsymbol{\hat{\alpha}}_{a}^{k}$, with the previous step, rather than using the current optimal solution to the decentralized DC-OPF. From the updated capacity price in \eqref{eq.PriceUpdate}, both sending and receiving areas use the same capacity price in their DC-OPF problem at each iteration $k$. Moreover, the capacity price is updated with a constant step size ($\beta$), while the inertial updates assume a diminishing step size. The convergence analysis is based on the exploitation of this fast update of the capacity prices versus the slow update of intertie flows and angles, and locational marginal prices. The update mechanism is complemented by an incentive mechanism (a money transfer and participation fee) in Step 4, to be analyzed in Section \ref{sec.inc} after discussing the convergence of the updates (Steps 1-3).

% This inertial update is essential for convergence of the iterative mechanism, because it means the capacity price is updated faster, and thus can catch up 

% \blue{From the updated capacity price in \eqref{eq.PriceUpdate}, both sending and receiving areas use the same capacity price in their DC-OPF problem at each iteration $k$. 
% Thus, considering that $\frac{\mu^{\ast}_{(i,j)}}{2}\mbox{sign} ({T}_{a,(i,j)}^*)=\bar{\eta}^{\ast}_{a,(i,j)}-\bar{\kappa}^{\ast}_{a,(i,j)}$, the consensus condition in Remark \ref{Rem.Consensus} is naturally met.}

\begin{algorithm}[h]
\SetAlgoLined
\SetKwInOut{Input}{Require}\SetKwInOut{Output}{Output}
\Input{Initialize $\mathbf{x}^0_a=(\mathbf{T}^0_a,\bm{\theta}^0_a,\bm{\alpha}^0_a)$, and $\mu^0_{(i,j)}$ for all areas and tielines.}
\Input{Maximum number of iterations $T+1$. Set $k=1$}
 \While{$k\leq T+1$}{
 %$|T^{k-1}_{a,(i,j)}+T^{k-1}_{a',(i,j)}|> \epsilon_1$ OR $|\mu^{k-1}_{i,j}-\mu^{k-2}_{i,j}|> \epsilon_2$
\begin{itemize}

\item Run DC-OPF problems of all areas and obtain $\hat{\mathbf{x}}^k_a=(\hat{\mathbf{T}}^k_a,\hat{\bm{\theta}}^k_a,\hat{\bm{\alpha}}^k_a)$.

\item Calculate $\mathbf{x}^k_a=(1-\rho_k)\mathbf{x}^{k-1}_a+\rho_k\hat{\mathbf{x}}^{k}_a$ for all areas.

\item Calculate the capacity prices $\mu _{(i,j)}^{k}$ \eqref{eq.PriceUpdate} for all tielines. %$\mu _{(i,j)}^{k}=\max \{\mu _{(i,j)}^{k-1}+\beta (\frac{\left\vert
% T_{a,(i,j)}^{k}\right\vert +\left\vert T_{a^{\prime },(i,j)}^{k}\right\vert 
% }{2}-\bar{T}_{a,(i,j)}),0\}$ for all tielines. 

\item Exchange $\mathbf{x}^k_a$ between neighboring areas. 
\item Set $k=k+1$.
\end{itemize}  
 }
 \caption{Iterative Tieline Capacity Allocation}
\label{al.Iterative}
\end{algorithm}

% \begin{figure}[ht]
% \centering
% \begin{tikzpicture}[node distance=5cm]

% \node (area1) [area] {AREA $A$ \[{\bf \hat{x}}_{a}^{k}=(\mathbf{\hat{T}}^{k}_{a},{\hat{\boldsymbol {\theta}}}^{k}_{a},\boldsymbol{\hat{\alpha}}^{k}_{a})\] 
% \begin{equation*}
% \mathbf{x}_{a}^{k} =  (1-\rho _{k})\mathbf{x}_{a}^{k-1}+\rho _{k} \hat{\mathbf{x}}_{a}^{k}
% \end{equation*}
% \begin{equation*} 
% \mu _{(i,j)}^{k}=\max \{\mu _{(i,j)}^{k-1}+\beta (\frac{\left\vert
% T_{a,(i,j)}^{k}\right\vert +\left\vert T_{a^{\prime },(i,j)}^{k}\right\vert 
% }{2}-\bar{T}_{a,(i,j)}),0\}
% \end{equation*}};
% \node (area2) [area, below of=area1] {AREA $A'$ \[\boldsymbol{\hat{x}}_{a'}^{k}=(\mathbf{\hat{T}}^{k}_{a'},\boldsymbol{\hat{\theta}}^{k}_{a'},\boldsymbol{\hat{\alpha}}^{k}_{a'})\] 
% \begin{equation*}
% \mathbf{x}_{a'}^{k} = (1-\rho _{k})\mathbf{x}_{a'}^{k-1}+\rho _{k} \hat{\mathbf{x}}_{a'}^{k}
% \end{equation*}
% \begin{equation*} 
% \mu _{(i,j)}^{k}=\max \{\mu _{(i,j)}^{k-1}+\beta (\frac{\left\vert
% T_{a,(i,j)}^{k}\right\vert +\left\vert T_{a^{\prime },(i,j)}^{k}\right\vert 
% }{2}-\bar{T}_{a,(i,j)}),0\}
% \end{equation*}};

% \draw [arrow] (area1) --  node[anchor=east] {$\mathbf{x}^{k-1}_a$} (area2);
% \draw [arrow] (area2) --node[anchor=west] {$\mathbf{x}^{k-1}_{a'}$}  (area1);
% \end{tikzpicture}

% \caption{\small Data flow for a two-area system at iteration $k.$ }
%  \label{fig.DataFlow}
% \end{figure}

\section{Convergence Analysis} \label{sec.convergence}

We make the following standing assumption on the feasibility of each area in the absence of tielines. 

\parindent 0cm
\medskip

% \textbf{Assumption (Maximum Intertie Capacity Price)}:
% {\em For each intertie $%
% (i,j)$, there exists a price $\bar{\mu}_{(i,j)}>0$ such that for $\mu
% _{(i,j)}\geq \bar{\mu}_{(i,j)}$ the optimal solution to \eqref{eq.Obj2} -  \eqref{eq.slack2} will have no flow along the intertie, i.e. $T_{a,(i,j)}=0$.}
% \red{CE: this assumption should hold whenever each area has enough generation capacity within. May be we can state this assumption as follows:
% There exists a feasible solution to \eqref{eq.Obj2}-\eqref{eq.slack2} for area $a$ when $T_{a,(i,j)}=0$ for all $(i,j) \in \mathcal{T}_a$.}

\textbf{Assumption (Area feasibility)}: There exists a feasible solution to the decentralized DC-OPF problem \eqref{eq.Obj2}-\eqref{eq.slack2} for area $a\in \mathcal A$ when $T_{a,(i,j)}=0$ for all $(i,j) \in \mathcal{T}_a$.
\medskip

The above assumption states that each market can meet the demand in its own area without relying on the intertie power flows. A manifestation of this assumption on coupled markets is the existence of a maximum intertie capacity price  $\bar{\mu}_{(i,j)}>0$ for $(i,j)\in \mathcal{T}_a$ such that for $\mu _{(i,j)}\geq \bar{\mu}_{(i,j)}$ the optimal solution to \eqref{eq.Obj2} -  \eqref{eq.slack2} will have no flow along the intertie, i.e. $T_{a,(i,j)}=0$---see Proposition \ref{prop_assumption} in the Appendix. 

In the first part of the convergence analysis, we show that the capacity multipliers converge.

\begin{lemma}\label{lem_capacity_price}
For every intertie $(i,j)$ it holds that
 $\mu _{(i,j)}^{k}\rightarrow \mu _{(i,j)}^{\infty }\geq 0$
and 
\begin{equation}\label{eq.compslack}
\mu _{(i,j)}^{k}(\frac{\left\vert T_{a,(i,j)}^{k}\right\vert +\left\vert
T_{a^{\prime },(i,j)}^{k}\right\vert }{2}-\bar{T}_{a,(i,j)})\rightarrow 0
\end{equation}
\end{lemma}
\begin{proof}
See appendix.
\end{proof}

The above result implies that the asymptotic capacity price $\mu_{(i,j)}^\infty$ is dual feasible and satisfies a complementary slackness condition \eqref{eq.compslack}.

\subsection{Convex Reporting Strategies}

% In other words, i.e. the limit outcome corresponds to the solution of the joint (centralized) DC-OPF for all areas with the reported supply and demand functions from each individual market clearing.}

We consider the possibility of strategic behavior by market participants in the form of  reporting strategies that are consistent with {\em some} choice of strictly convex generation cost functions. Specifically, in the course of the execution of the trading mechanism, we consider the possibility that area $a$ reports values of the form 
$
\boldsymbol{\hat{x}}_{a}^{D,k}=(\mathbf{\hat{T}}^{D,k}_{a},%
\boldsymbol{\hat{\theta}}^{D,k}_{a},\boldsymbol{\hat{\alpha}}^{D,k}_{a})
$ 
that corresponds to the optimal solution of the problem:
\begin{align} 
\min & ~~~\sum_{g\in \mathcal{G}_a} C^{D}_{a,g}(P_{a,g})+\sum_{(i,j)\in \mathcal{T}_a} \big[\frac{\mu_{(i,j)}^{k-1}}{2}|T_{a,(i,j)}|-\alpha^{k-1}_{a',j} T_{a,{(i,j)}}\big],
\label{deviation} \\
& ~~~\mbox{s.t. (\ref{eq.Bal2}), (\ref{eq.Plim2}), (\ref{eq.Flim2}),   (\ref{eq.slack2})} \nonumber \\
& ~~~ T_{a,(i,j)}= \frac{\theta_{a,i}-\theta_{a',j}^{k-1}}{\bar{x}_{a,(i,j)}}, ~~\forall  (i,j) \in \mathcal{T}_a,
\nonumber
\end{align}

where $C_{a,g}^{D}(\cdot)$ are strictly convex differentiable functions for each $g \in \mathcal {G}_a$. We shall refer to this form of reporting as {\em convex}. 

\medskip
{\bf Definition 1} {\em (Convex reporting strategy) Area $a \in \mathcal{A}$ follows a convex reporting strategy if and only for every iteration $k>0$, such area reports the values 
$
\boldsymbol{\hat{x}}_{a}^{D,k}=
(\mathbf{\hat{T}}^{D,k}_{a},
\boldsymbol{\hat{\theta}}^{D,k}_{a},
\boldsymbol{\hat{\alpha}}^{D,k}_{a})
$ 
corresponding to the solution to (\ref{deviation}) with $C_{a,g}^{D}(\cdot ), g \in \mathcal{G}_a$ strictly convex.}

%A convex reporting strategy by an area is not truthful when $C_{a,g}^{D}(\cdot) \neq C_{a,g}(\cdot)$. We define the truthful reporting strategy formally as follows. 

\medskip
Our main convergence result shows that the sequence of outcomes of the proposed market coupling mechanism converges to the optimal allocation and pricing of intertie capacity {\em given} the reported supply and demand functions from each individual market clearing. 

% We now state and prove the main convergence result.

\begin{theorem}\label{thm_convergence}
Assuming each area $a\in \mathcal A$ follows a convex reporting strategy with $C_{a,g}^{D}(\cdot )=C_{a,g}(\cdot ), g \in \mathcal{G}_a$, the sequence 
$\boldsymbol{x}_{a}^{k} \triangleq (\mathbf{T}_{a}^{k},\boldsymbol{\theta}_{a}^{k},\boldsymbol{\alpha}_{a}^{k})$, $a \in \mathcal{A}$ generated by Algorithm \ref{al.Iterative} 
converges to intertie flows, angles and locational marginal prices in the centralized DC-OPF problem.
\end{theorem}
\begin{proof}
The proof has two parts. In the first part, we show that under diminishing step size updates, the sequence $\bbx_a^k$ is Cauchy, and thus $\bbx_a^k$ converges to $\bbx_a^\infty$. In the second part, we use this convergence along with the convergence of capacity prices (Lemma \ref{lem_capacity_price}) to show that the KKT conditions of the centralized and decentralized DC-OPF are equivalent, concluding that $\bbx_a^\infty$ and $\mu^\infty_{(i,j)}$ correspond to the centralized optimal solution.

By definition of the update rule, it holds that
\begin{align*}
&\boldsymbol{x}_{a}^{k+1}-\boldsymbol{x}_{a}^{k} \\&=(1-\rho _{k+1})%
\boldsymbol{x}_{a}^{k}+\rho _{k+1}\boldsymbol{\hat{x}}_{a}^{k+1}-(1-\rho
_{k})\boldsymbol{x}_{a}^{k-1}-\rho _{k}\boldsymbol{\hat{x}}_{a}^{k} \\
&=(1-\rho _{k})(\boldsymbol{x}_{a}^{k}-\boldsymbol{x}_{a}^{k-1})+\rho _{k}(%
\boldsymbol{\hat{x}}_{a}^{k+1}-\boldsymbol{\hat{x}}_{a}^{k}) +(\rho
_{k+1}-\rho _{k})(\boldsymbol{\hat{x}}_{a}^{k+1}-\boldsymbol{x}_{a}^{k})
\end{align*}

Since the optimal solution and Lagrange multipliers for each area $a \in \cal{A}$ are Lipschitz
continuous (see e.g. \cite{Robinson}), it follows that%
\begin{align*}
\left\Vert \boldsymbol{\hat{x}}_{a}^{k+1}-\boldsymbol{\hat{x}}_{a}^{k}\right\Vert &\leq \sum\limits_{a^{\prime }\in \mathcal{A}%
(a)}L_a\left\Vert \boldsymbol{x}_{a^{\prime }}^{k+1}-\boldsymbol{x}_{a^{\prime
}}^{k}\right\Vert + L_a\|\boldsymbol{\mu}_{a}^{k}-\boldsymbol{\mu}_{a}^{k-1}\| \\&=\mathcal{O}(\rho _{k})+L_a\|\boldsymbol{\mu}_{a}^{k}-\boldsymbol{\mu}_{a}^{k-1}\|
\end{align*}
for some $L_a>0$. Thus for finite $T<\infty $ it follows that%
\begin{align*}
\left\Vert \boldsymbol{x}_{a}^{k+T}-\boldsymbol{x}_{a}^{k}\right\Vert  &\leq
\sum\limits_{\ell =1}^{T}\Bigg[ \rho _{k+\ell }\left\Vert \boldsymbol{x}%
_{a}^{k+\ell }-\boldsymbol{x}_{a}^{k+\ell -1}\right\Vert 
+\rho _{k+\ell}\left\Vert \boldsymbol{\hat{x}}_{a}^{k+\ell +1}-\boldsymbol{\hat{x}}%
_{a}^{k+\ell }\right\Vert \nonumber  +\left\vert \rho _{k+\ell +1}-\rho _{k+\ell}\right\vert \left\Vert \boldsymbol{\hat{x}}_{a,i}^{k+\ell +1}-\boldsymbol{x}_{a,i}^{k+\ell }\right\Vert \Bigg]  
\\ & \leq  \sum\limits_{\ell =1}^{T}\rho _{k+\ell }\mathcal{O}(\rho _{k+\ell
})+\sum\limits_{\ell =1}^{T}\sum\limits_{(i,j)\in \mathcal{T}_{a}}\rho
_{k+\ell }(\mu _{(i,j)}^{k+\ell}-\mu _{(i,j)}^{k+\ell-1})+ \mathcal{O}(1)\sum\limits_{\ell =1}^{T}|\rho _{k+\ell +1}-\rho _{k+\ell }|.
\end{align*}%
Thus the sequence $\{\boldsymbol{x}_{a}^{k}:k>0\}$ is Cauchy and converges
to say $\boldsymbol{x}_{a}^{\infty }$. Since $\boldsymbol{x}_{a'}^{k}\rightarrow \boldsymbol{x}_{a'}^{\infty }$ for all adjacent areas $a' \in \mathcal{A}$, by continuity of optimal solutions and
Lagrange multipliers it follows that the sequence $\{\boldsymbol{\hat{x}}%
_{a}^{k}:k>0\}$ has a limit. By construction, $\boldsymbol{x}%
_{a}^{k}$ is a weighted average of $\{\boldsymbol{\hat{x}}%
_{a}^{\ell}\}_{\ell =0}^{k}$, hence $\|\boldsymbol{\hat{x}}%
_{a}^{k}-\boldsymbol{x}%
_{a}^{k}\|\rightarrow 0$ and $\boldsymbol{\hat{x}}%
_{a}^{k} \rightarrow \boldsymbol{x}%
_{a}^{\infty}$.

In the second part of the proof we show that the iterative mechanism satisfies the KKT conditions of the centralized OPF problem. 

We begin by showing that $\{\boldsymbol{x}^\infty\}_{a\in\ccalA}$ is a feasible solution of  the centralized OPF problem. Note that a decentralized DC-OPF solution that is feasible for all areas satisfies all the conditions of the centralized DC-OPF except \eqref{eq.Tielim1}. The convergence of the sequence $\bbx^k$ to $\bbx^\infty$ implies that $T_{a,(i,j)}^\infty = -T_{a',(i,j)}^\infty$ for $(a,a')\in \ccalT$. We consider two cases of the capacity price: $\mu_{(i,j)}=0$ and $\mu_{(i,j)}>0$ to argue that inter-tie capacity limits are obeyed. 
If $\mu_{(i,j)}^\infty=0$, then it must be that $T_{a,(i,j)}^\infty \leq \bar T_{a,(i,j)}$ because of the capacity price updates in \eqref{eq.PriceUpdate}. If $\mu_{(i,j)}^\infty>0$, then $T_{a,(i,j)}^\infty=\bar T_{a,(i,j)}$ by the condition in \eqref{eq.compslack}. Thus, the intertie capacity limits are satisfied by $\boldsymbol{x}_a^\infty$. The dual feasibility and the complementary slackness conditions hold by Lemma \ref{lem_capacity_price}. 

Next we show that the first order stationarity conditions of the centralized problem are satisfied by $\mathbf{T}_{a}^{\infty}$ and $\boldsymbol{\theta}_{a}^{\infty}$. In the limit as $k \rightarrow \infty$, let $\mathcal{L}_{A}(\boldsymbol{x}_A,\boldsymbol{x}_{B}^\infty,\alpha_{B}^{\infty},\mu^\infty)$ and $\mathcal{L}_{B}(\boldsymbol{x}_{B},\boldsymbol{x}_{A}^\infty,\alpha_{A}^{\infty},\mu^\infty)$ denote the Lagrangian of the decentralized DC-OPF problem for areas $A$ and $B$, respectively. The terms  $\boldsymbol{x}_{B}^\infty,\alpha_{B}^{\infty}$, and $\mu^\infty$ in $\mathcal{L}_{A}(\cdot)$ represent the terms in the objective \eqref{eq.Obj2} and tie-flow constraint \eqref{eq.TieFlow2} of Area A coming from the decoupling of constraints in the centralized DC-OPF. We note that $\frac{\partial \mathcal{L}_A(\cdot)}{\partial\mathbf{P}_A}=\frac{\partial \mathcal{L}(\cdot)}{\partial\mathbf{P}_A}$ and $\frac{\partial \mathcal{L}_A(\cdot)}{\partial \theta_{A,i}}=\frac{\partial \mathcal{L}(\cdot)}{\partial\theta_{A,i}}$ $\forall (i,j) \in \mathcal{F}_a$ where $\mathcal{L}(\cdot)$ is the Lagrangian of the centralized DC-OPF. 

Next we consider the first order stationarity conditions with respect to $\theta _{A,i}$, $\theta _{B,j}$ 
and $T_{A,(i,j)}$ for tieline $(i,j)\in \mathcal{T}_{A}$ in the decentralized DC-OPF problem and show their equivalence to the corresponding first order conditions of the centralized DC-OPF in  \eqref{eq.centralizedangle_1}-\eqref{eq_centralized_intertie}.
%
% At iteration $k>0$, the first order
% conditions with respect to $T_{a,(i,j)}$ are 
% \begin{equation}
% \hat{\alpha}_{a,i}^{k}+\hat{\xi}_{a,(i,j)}^{k}=\alpha_{a^{\prime },j}^{k-1}-\frac{%
% {\mu _{(i,j)}^{k-1}}}{2}\text{sign}(\hat{T}_{a,(i,j)}^{k})
% \label{optimality_a}
% \end{equation}%
% Similarly, for area $a^{\prime }$ the first order conditions with respect to $T_{a^{\prime },(i,j)}$ are%
% \begin{equation}
% \hat{\alpha}_{a^{\prime },j}^{k}+\hat{\xi}_{a^{\prime },(i,j)}^{k}=\alpha
% _{a,i}^{k-1}-\frac{\mu _{(i,j)}^{k-1}}{2}\text{sign}(\hat{T}_{a^{\prime
% },(i,j)}^{k}) 
% \label{optimality_a_prime}
% \end{equation}%
%
The first order conditions with respect to $T_{a,(i,j)}$ for areas $a\in \{A,B\}$, i.e., $\frac{\partial \mathcal{L}_{A}(\boldsymbol{x}_A,\boldsymbol{x}_{B}^\infty,\alpha_{B}^{\infty},\mu^\infty)}{\partial T_{A,(i,j)}}$ and $\frac{\partial \mathcal{L}_{B}(\boldsymbol{x}_B,\boldsymbol{x}_{A}^\infty,\alpha_{A}^{\infty},\mu^\infty)}{\partial T_{B,(i,j)}}$, are respectively as follows 
\begin{eqnarray}
\hat \alpha_{A,i}+\hat{\xi}_{A,(i,j)} &=&\alpha _{B,j}^{\infty }-\frac{{\mu_{(i,j)}^{\infty }}}{2}\text{sign}(\hat T%
_{A,(i,j)}), \label{eq_dec_kkt_a}\\
\hat \alpha _{B,j}+\hat{\xi}_{B,(i,j)}
&=&\alpha _{A,i}^{\infty }-\frac{\mu_{(i,j)}^{\infty }}{2}\text{sign}(\hat T
_{B,(i,j)}). \label{eq_dec_kkt_b}
\end{eqnarray}%
From Lemma \ref{lemma_equivalence}, we know that if $\alpha _{B,j}^{\infty} = \alpha _{B,j}^{*}$ and $\frac{\mu _{(i,j)}^{\infty}}{2} \sign(\hat T_{A,(i,j)}) = \bar \eta_{A,(i,j)}^*- \bar \kappa_{A,(i,j)}^*$, then the centralized solution ($\hat \alpha _{A,i}=\alpha _{A,i}^*$ and $\hat{\xi}_{A,(i,j)}={\xi}_{A,(i,j)}^*$) satisfies the stationarity condition in \eqref{eq_dec_kkt_a} for the decentralized DC-OPF for area A. The same argument applies to the condition in \eqref{eq_dec_kkt_b} for area B. 
Consider the optimal dual variables ($\bar \eta_{A,(i,j)}^*$ and $\bar \kappa_{A,(i,j)}^*$) associated with the intertie capacity constraints \eqref{eq.Tielim1}. The intertie capacity price 
($\frac{\mu_{(i,j)}^\infty}{2} \sign(T_{A,(i,j)}^\infty)$) satisfies the identical complementary slackness conditions as $\bar \eta_{A,(i,j)}^*-\bar \kappa_{A,(i,j)}^*$. Same applies to the optimal dual variables of area B ($\bar \eta_{B,(i,j)}^*$ and $\bar \kappa_{B,(i,j)}^*$). 
Thus, we have $\frac{\partial \mathcal{L}_A(\boldsymbol{x}_A^*,\boldsymbol{x}_B^\infty,\alpha_B^\infty,\mu_{(i,j)}^\infty)}{\partial T_{A,(i,j)}}=0$ if $\alpha_B^\infty=\alpha_B^*$ and $\frac{\partial \mathcal{L}_B(\boldsymbol{x}_B^*,\boldsymbol{x}_A^\infty,\alpha_A^\infty,\mu_{(i,j)}^\infty)}{\partial T_{B,(i,j)}}=0$, if $\alpha_A^\infty=\alpha_A^*$. For these two conditions to hold at the same time, it must be that $\boldsymbol{x}_A^\infty=\boldsymbol{x}_A^*$ and $\boldsymbol{x}_B^\infty=\boldsymbol{x}_B^*$.
\end{proof}

The main steps in the proof of Theorem \ref{thm_convergence} involves showing convergence of the capacity prices (Lemma \ref{lem_capacity_price}) and the intertie flows, angles, and locational marginal prices ($\boldsymbol{x}_{a}^{k}$) to some finite value. In showing their convergence, we rely on the updates of $\boldsymbol{x}_{a}^{k}$ with diminishing step size in \eqref{eq_inertia_a}-\eqref{eq_inertia_c} and Lipschitz continuity of solutions to the decentralized DC-OPF problem. In the second part of the proof, we show that these convergence points have to correspond to the optimal solution of the centralized DC-OPF problem via the equivalence of KKT conditions between the decentralized and centralized DC-OPF problems. 

In the remainder of the paper we examine incentive compatibility issues that may be induced by market coupling.

\section{An Incentive Compatible Design for Market Coupling} \label{sec.inc}
In this section, we identify incentive transfers based on {\em marginal contribution to the coupling} of each area to ensure that market coupling does not introduce incentives for strategic manipulation of reported supply cost functions.
To focus on incentive compatibility issues that may exclusively arise from market coupling, we assume the reported supply and demand functions for the clearing of each market, i.e.  $C_{a,g}(\cdot ), g \in \mathcal{G}_a$  constitute a {\em Nash equilibrium} for internal market clearing. Since each area (or market) is only an intermediary, the net savings resulting from coupling are {\em fully distributed} back to market participants. Hence, to show that there is no incentive to manipulate the reported supply cost functions, it is sufficient to show that any manipulation  of reported information does not yield increased net savings at the {\em aggregate} market level as a result of the coupling of all markets.

\medskip
The marginal contribution of an area $a\in\mathcal{A}$ refers to the cost savings of other areas 
$a'\in \mathcal{A}\backslash a$ as a result of area $a$ participating in the coupling. 

To formalize this idea, let $(\boldsymbol{P}^{*}_{a},\boldsymbol{P}^{*}_{-a})$ denote the solution identified in Theorem 1 with all areas participating with the reported supply cost functions for internal market clearing and
\begin{equation*}
C^{*}_a:=\sum_{g\in \mathcal{G}_{a}} C_{a,g}(P^{*}_{a,g}) 
\end{equation*}
$C^{*}:=\sum_{a \in \mathcal{A}}C^{*}_a$ and $C^{*}_{-a}:=\sum_{a' \in \mathcal{A}\backslash a}C^{*}_{a'}$. Similarly, let $(\boldsymbol{P}_{a,0}, \boldsymbol{\tilde{P}}_{-a})$ denote the solution of the DC-OPF problem {\em without} area $a$'s participation (i.e. area $a$'s initial power flow and exchange with adjacent areas remains fixed). With 
$\tilde C_{a'}:=\sum_{g\in \mathcal{G}_{a'}} C_{a',g}(\tilde{P}_{a',g})$, the change in cost induced by area $a$'s participation to all other areas in the coupling can be expressed as 
$$
M_a := C^*_{-a} -\tilde{C}_{-a}
$$ where $\tilde{C}_{-a}:=\sum_{a'\in \mathcal{A} \backslash a} \tilde{C}_{a'}$.

\medskip
In what follows, we will describe an iterative design with incentive transfers that upon stopping (approximately) gives each area its marginal contribution to the coupling. Thus, upon stopping, area $a$ is expected to have a net cost reduction equal to:
\begin{align*}
    C^*_{a}- C_{a,0} + M_a &= C^*_{a} -C_{a,0} + C^*_{-a} -\tilde{C}_{-a}\\
&=
C^*_{a} + C^*_{-a}-(C_{a,0}+\tilde{C}_{-a}) \\
&= C^{\ast}-(C_{a,0}+\tilde{C}_{-a})  < 0
\end{align*}
where $C_{a,0}$ is the cost of area $a$ under the initial power flow. The inequality follows from the fact that the combination $(\boldsymbol{P}_{a,0},\boldsymbol{\tilde{P}}_{-a})$ of area $a$'s initial dispatch $\boldsymbol{P}_{a,0}$ and the resulting coupled dispatch $\boldsymbol{\tilde{P}}_{-a}$ is a feasible power dispatch for the interconnected power system with all areas while $\boldsymbol{P}^{\ast}_{a}$ is the optimal solution to the interconnected DC-OPF with all areas (Theorem 1) given internal market clearing conditions.
%An unsustainable incentive compatible design can ensure each area is truthful by dispensing money into the coupled market.
To ensure incentive transfers and cost savings are balanced, we introduce a minimum {\em participation} fee $R>0$ that is assessed to each market (area) participating in the proposed design for market coupling.

\subsection{Marginal Contribution to a Coupling}

We will use $\boldsymbol{\tilde{x}}^{k}=(\mathbf{\tilde{T}}^{k},\boldsymbol{\tilde{\theta}}^{k},\boldsymbol{\tilde{\alpha}}^{k})$ to denote the sequence of the flows between areas, angles and nodal prices for relevant nodes in adjacent areas at iteration $k$ generated by the updates in Algorithm \ref{al.Iterative} without area $a$'s participation in the coupling. The cost of area $a'$ associated with the updates $\boldsymbol{\tilde{x}}^{k}$ is denoted using $\tilde C_{a',k}:=\sum_{g\in \mathcal{G}_{a'}} C_{a',g}(\tilde{P}_{a',g}^{k})$. Similarly we use $C_{a',k}$ to denote the cost of area $a'$ associated with the sequence $\boldsymbol{\hat{x}}^{k}$ when area $a$ is participating in the coupling. The incremental change in cost of area $a'$ from iteration $k$ to iteration $k+1$ with and without area $a$'s participation are respectively given as $\Delta \tilde C_{a',k}:=\tilde C_{a',k+1}-\tilde C_{a',k}$ and $\Delta C_{a',k}:= C_{a',k+1}- C_{a',k}$. The difference between the incremental changes in costs $\Delta \tilde C_{a',k}-\Delta C_{a',k}$ represents the incremental cost saving for area $a' \in\mathcal{A} \setminus a$ caused by area $a$'s participation in the iterative mechanism. Formally, we define the marginal contribution of area $a$ at step $k$ as the change in cost to all areas $a'\in \ccalA\setminus\{a\}$ relative to the change in cost when \emph{area a is not part of the coupling},
\begin{equation} \label{eq_marginal}
    \Delta M_{a,k} := \sum_{a'\in \ccalA \setminus \{a\}} \Delta C_{a',k}-\Delta \tilde C_{a',k}.
\end{equation}
Using the definitions above, we define the incentive transfer for area $a$ as follows:
\begin{equation} \label{eq_incentive}
\Delta \pi _{a,k}=
\Delta M_{a,k}+r_{a,k+1}-r_{a,k}
\end{equation}
where
\begin{equation} \label{eq.rak}
r_{a,k}:=\hspace{-8pt}\sum_{(i,j)\in \mathcal{T}_{a}}{\alpha _{a^{\prime },j}^{k-1}}\hat{T}_{a,(i,j)}^{k}-\hspace{-8pt}\sum_{(i,j)\in \mathcal{T}_{a}}\frac{\mu _{(i,j)}^{k-1}}{2}%
(\left\vert \hat{T}_{a,(i,j)}^{k}\right\vert -\bar{T}_{a,(i,j)}).
\end{equation}
The incentive transfers give area $a$ its marginal contribution to the coupling and offsets the net change value of intertie trades for area $a$, i.e. $r_{a,k+1}-r_{a,k}$.

\begin{figure}[h] 
	\centering
\vspace{-5pt}
	\includegraphics[width=0.5\linewidth]{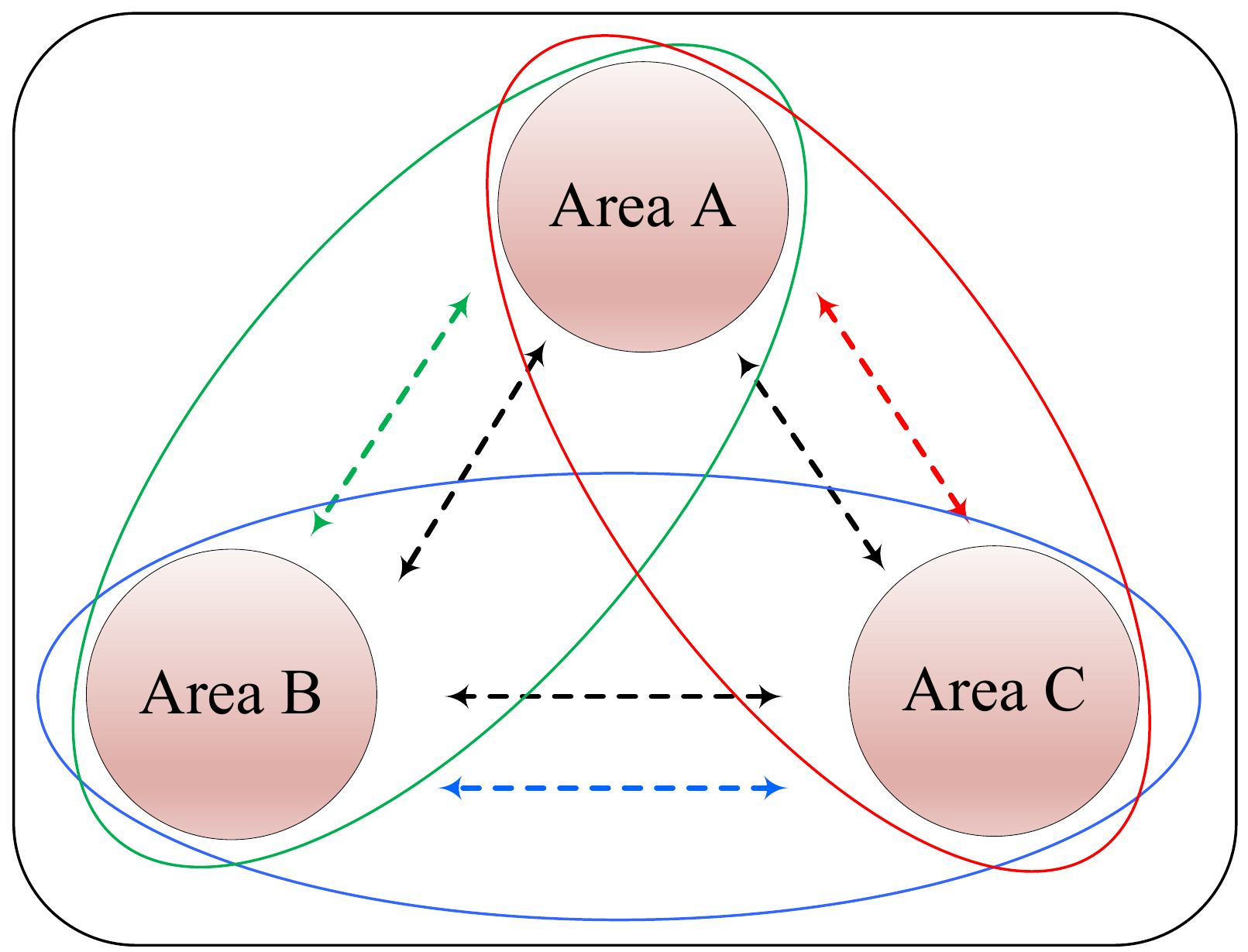}
    \vspace{-5pt}
	\caption{\small Marginal contribution to a coupling. The dashed and solid lines respectively represent the data exchange and coupling territory.
The black lines refer to the case where all areas participate in market coupling while the blue, red, and green lines are respectively attributed to the cases where either area A, B, or C is excluded from coupling.}
       \label{fig.MC}
\end{figure}

\medskip
\begin{remark}
 The computation of $\Delta \pi_{a,k}$ for each $a \in \mathcal{A}$
requires a total of $|\mathcal{A}+1|$ parallel executions of Algorithm 1: one for the complete coupling of all areas and $|\mathcal{A}|$ for different coupling excluding each one of the areas.
% Computation of $\Delta \tilde{C}_{a^{\prime },k}$ in the incentive transfers  requires each area to run $|\mathcal{A}|$ updates in parallel considering all hypothetical scenarios with an area excluded from the coupling. 
Figure \ref{fig.MC} shows all possible subsets of coupling with one area excluded from the coupling of three areas $\ccalA=\{A,B,C\}$. For instance, Areas B and C need to exchange information and bids according to Algorithm \ref{al.Iterative} as if area A is not part of the coupling---see blue lines in the figure. This will generate sequence of nodal prices, angles and intertie flows, denoted by $\boldsymbol{\tilde{x}}^{k}=\{(\mathbf{\tilde{T}}^{k}_a,\boldsymbol{\tilde{\theta}}^{k}_a,\boldsymbol{\tilde{\alpha}}^{k}_a)\}_{a\in\ccalA}$,  different from the sequence $\boldsymbol{\hat{x}}^{k}=\{(\mathbf{\hat{T}}^{k}_{a},{\hat{\boldsymbol {\theta}}}^{k}_{a},\boldsymbol{\hat{\alpha}}^{k}_{a})\}_{a\in \ccalA}$ generated with all areas included in the coupling. %Indeed, the system needs to generate $|\ccalA|$ sequences, each sequence representing a single area excluded from the coupling in order to compute the incentive transfer and participation fees. 
In practice, the number of markets in the coupling would be small for these hypothetical updates to be computationally demanding.
\end{remark}
\medskip

\subsection{Estimation of Marginal Contribution to Coupling}
The incentive transfers defined \eqref{eq_incentive} can not be readily computed since the change in costs $\Delta C_{a,k}$ is privately known by area $a$'s operator. However, as we shall show below a fairly accurate estimate of an area's marginal contribution to a coupling $\widehat{M_{a}}$ can be obtained using only reported information (i.e. desired transfer flows, angles and locational marginal prices) for each area $a \in \mathcal{A}$. 

Let $V_{a,k}$ denote the optimal value of area $a$'s DC-OPF problem %
\eqref{eq.Obj2}-\eqref{eq.slack2} at iteration $k$, 
\begin{equation}
V_{a,k}=C_{a,k}-r_{a,k}.  \label{eq_optimal_k}
\end{equation}
We use the first order term Taylor expansion of the optimal value to define the incentive payments at each step,
\begin{equation}
V_{a,k+1}=V_{a,k}+\Delta V_{a,k}+\mathcal{O}(\beta ^{2})  \label{Taylor}
\end{equation}
From \eqref{eq_optimal_k}, it follows that $\Delta C_{a,k}=\Delta V_{a,k} + r_{a,k+1}- r_{a,k}+\mathcal{O}(\beta^2)$. As we show below, the value $\Delta V_{a,k}$ can be estimated based upon reported information. This allows us to estimate the value of $\Delta C_{a,k}$.

 \begin{lemma} \label{lemma_deltaV}
A first order approximation to the change in optimal value for area $a$ is given by: 
\begin{align}
\label{change_V}
\Delta &V_{a,k}=-\hspace{-1pt}\sum\limits_{(i,j)\in \mathcal{T}_{a}}\Big[ \hat{T}%
_{a,(i,j)}^{k}({\alpha _{a^{\prime },j}^{k}-\alpha _{a^{\prime },j}^{k-1})}%
\hspace{6pt}-\frac{1}{2}\left\vert \hat{T}_{a,(i,j)}^{k}\right\vert (\mu _{(i,j)}^{k}-\mu
_{(i,j)}^{k-1})-\frac{\hat{\xi}_{a,(i.j)}^{k}}{\bar{x}_{a,(i.j)}}(\theta
_{a^{\prime },j}^{k}-\theta _{a^{\prime },j}^{k-1})\Big]
\end{align}
with:
\begin{equation}
\label{xi_opt}
\hat{\xi}_{a,(i,j)}^{k}=\alpha _{a^{\prime },i}^{k-1}-\frac{\mu
_{(i,j)}^{k-1}}{2}\text{sign}(\hat{T}_{a,(i,j)}^{k})-\hat{\alpha}_{a,i}^{k}.
\end{equation}
\end{lemma}
\begin{proof} 
Change in optimal value is given by 
\begin{align*}
\Delta V_{a,k}&\triangleq \sum_{(i,j)\in \mathcal{T}_{a}} \Big[ \frac{\partial V_{a,k}}{\partial {\alpha _{a^{\prime },j}}}%
(\alpha_{a^{\prime},j}^{k}-\alpha _{a^{\prime },j}^{k-1}) %
\hspace{6pt}+\frac{\partial V_{a,k}}{\partial%
\mu _{(i,j)}} (\mu _{(i,j)}^{k}-\mu _{(i,j)}^{k-1}) +\frac{\partial V_{a,k}}{%
\partial \theta _{a^{\prime },j}}(\theta _{a^{\prime },j}^{k}-\theta%
_{a^{\prime },j}^{k-1}) \Big].
\end{align*}
By the envelope theorem, we have
 \begin{equation*}
\begin{array}{rrr}
\frac{\partial V_{a,k}}{\partial {\alpha _{a^{\prime },j}}}=-\hat{T}%
_{a,(i,j)}^{k}, &   \frac{\partial V_{a,k}}{\partial \mu _{(i,j)}}=
\frac{1}{2}\left\vert \hat{T}_{a,(i,j)}^{k}\right\vert, & \frac{\partial V_{a,k}}{\partial \theta _{a^{\prime},j}}
=\frac{\hat{\xi}_{a,(i.j)}^{k}}{\bar{x}_{a,(i.j)}}.
\end{array}
 \end{equation*}
Equation \eqref{xi_opt} is the first order optimality condition for flow accross intertie $(i,j)$ for market $a$.
\end{proof}

Using \eqref{eq_optimal_k} and \eqref{Taylor}, the change in cost for area $a' \in \mathcal{A}$ can be written as:
\begin{eqnarray*}
\Delta C_{a',k} %&\:= &\sum_{g\in \mathcal{G}_{a}}C_{a,g}(\hat{P}_{a,g}^{k+1})-\sum_{g\in \mathcal{G}_{a}}C_{a,g}(\hat{P}_{a,g}^{k}) \\
&=&V_{a',k+1}-V_{a',k}+r_{a',k+1}-r_{a',k} \notag \\
&=&\Delta V_{a',k}+r_{a',k+1}-r_{a',k}+\mathcal{O}(\beta ^{2}) \label{delta_C}.
\end{eqnarray*}%
We define the estimated changes in cost for area $a'$ by ignoring the $\mathcal{O}(\beta ^{2})$ terms, i.e.,
\begin{equation} \label{estimate_delta_C}
    \widehat{\Delta C}_{a',k} = \Delta V_{a',k}+r_{a',k+1}-r_{a',k}.
\end{equation}
Similarly, we obtain the following expression for the scenario with area $a$ excluded from the coupling, i.e.,
\begin{equation*}\label{estimate_delta_C_prime}
\Delta\tilde C_{a',k}=\tilde{\Delta V_{a',k}}+\tilde r_{a',k+1}-\tilde r_{a',k} +\mathcal{O}(\beta ^{2})
\end{equation*}
and an estimate of $\Delta \tilde{C}_{a',k}$:
\begin{equation}\label{estimate_delta_C_prime}
\widehat{\Delta \tilde{C}}_{a',k} =\tilde{\Delta V_{a',k}}+\tilde r_{a',k+1}-\tilde r_{a',k} 
\end{equation}
Thus, using \eqref{estimate_delta_C} and \eqref{estimate_delta_C_prime}, the estimate of an area's marginal contribution to a coupling can be readily computed using publicly available information, i.e., locational marginal prices, and intertie flows and angles as:
 
\[\Delta \widehat{M_{a,k}}=
\sum_{a'\in \ccalA \setminus \{a\}} \widehat{\Delta C}_{a',k}-\widehat{\Delta \tilde C}_{a',k}
\]
The incentive transfer for area $a$ is defined as:
\begin{equation} \label{incentive_hat}
\widehat{\Delta \pi}_{a,k}=\Delta \widehat{M_{a,k}} + r_{a,k+1}-r_{a,k}.
\end{equation}

\begin{remark}
The incentive transfer in \eqref{incentive_hat} is similar to the Clarke pivot rule in Vickrey-Clarke-Groves (VCG) mechanism \cite{Clarke},\cite{Vickrey}, \cite{Groves} as it applies to electricity markets \cite{NaLi} and \cite{Karaca}. However, unlike the VCG mechanism (which requires agents to report high dimensional data and a centralized solution of social surplus maximization), in the proposed market coupling reported information (i.e. energy flows and locational marginal price at interties) is elicited in an iterative {\em pointwise} manner and only a finite number of queries are required to approximate in a decentralized fashion the solution of the optimal allocation and pricing of intertie capacity.  
\end{remark}

\subsection{Net Benefit from Coupling}
Our first result identifies a sufficient condition on the participation fee $R>0$ to ensure all areas stand to obtain (approximately) a surplus (cost reduction) from participation.

%%%%%%%%%%%%%%%%%%%%%%%%%%%%%%%%%%%%%%%%%%%%%%%%%%%%%%%%%%%
%%%%%%%% BEGIN PROPOSITION %%%%%%%%%%%%%%%%%%%%%%%%%%%%%%%%%%%%%%
%%%%%%%%%%%%%%%%%%%%%%%%%%%%%%%%%%%%%%%%%%%%%%%%%%%%%%%%%%%

\begin{proposition} \label{prop_rational}
Consider the mechanism defined by the approximate incentive payments $\widehat{\Delta \pi}_{a,k}$ in \eqref{incentive_hat} and assume the participation fee satisfies the following condition:
\begin{equation}
    R \leq \min_{a \in \mathcal{A}}\{C_{a,0}+\tilde{C}_{-a}-C^{\ast}\}.
    \label{eq_participation_fee_prime}
\end{equation}
 For each area $a \in \mathcal{A}$, participation in the mechanism yields a net cost reduction bounded by $\delta _{a,T}(\beta)\geq0$ i.e.,
\begin{equation*}
V_{a,T+1}-V_{a,0}+\sum_{k=0}^{T}\widehat{\Delta \pi}_{a,k}+R \leq \delta_{a,T}(\beta)
\end{equation*} where T denotes the iteration at which the algorithm stops. Further, $\delta_{a,T}(\beta)\rightarrow 0$ as $\beta \rightarrow 0^+, T \rightarrow \infty$.

\end{proposition} 
%%%%%%%%%%%%%%%%%%%%%%%%%%%%%%%%%%%%%%%%%%%%%%%%%%%%%%%%%%%
%%%%%%%% END PROPOSITION %%%%%%%%%%%%%%%%%%%%%%%%%%%%%%%%%%%%%%
%%%%%%%%%%%%%%%%%%%%%%%%%%%%%%%%%%%%%%%%%%%%%%%%%%%%%%%%%%%
\begin{proof}
At the $k$-th iteration, the change in optimal value inclusive of incentive
transfer for area $a$ is%
\begin{align*}
V_{a,k+1}-V_{a,k}+\widehat{\Delta \pi} _{a,k} &= C_{a,k+1}-C_{a,k} +\hspace{-8pt}\sum_{a^{\prime
}\in \left. \mathcal{A}\right\backslash \{a\}}[\widehat{\Delta C}_{a^{\prime
},k}-\widehat{\Delta\tilde{C}}_{a^{\prime },k}], 
\end{align*}
where the equality follows by the cancellation of the term $r_{a,k+1}-r_{a,k}$  when we substitute in \eqref{eq_optimal_k} for $V_{a,k}$ and \eqref{incentive_hat} for $\widehat{\Delta \pi}_{a,k}$.
\begin{align}
   V_{a,k+1}-V_{a,k}&+\widehat{\Delta \pi} _{a,k} =C_{a,k+1}-C_{a,k}+C_{-a,k+1}-C_{-a,k}\hspace{2pt} -(\tilde{C}_{-a,k+1}-\tilde{C}%
_{-a,k})+2\left\vert \mathcal{A}-1\right\vert \mathcal{O}(\beta ^{2}).  \notag
\end{align}
where
\begin{equation*}
\begin{array}{rrr}
C_{a,k}=\sum\limits_{g\in \mathcal{G}_{a}}C_{g}(\hat{P}_{a,g}^{k}) &  & 
C_{-a,k}=\sum\limits_{a^{\prime }\in \left. \mathcal{A}\right\backslash
\{a\}}C_{a^{\prime },k}%
\end{array}%
\end{equation*}%
Upon stopping after $T$-iterations the total change in cost inclusive of
incentive transfer and participation fee for area $a$ is
\begin{align*}
\sum_{k=0}^{T}[V_{a,k+1}-V_{a,k}+&\widehat{\Delta \pi} _{a,k}] + R \\
&= C_{a,T+1}-C_{a,0}+C_{-a,T+1}-C_{-a,0}-\tilde{C}_{-a,T+1} +\tilde{C}_{-a,0}+ 2T\left\vert \mathcal{A}-1\right\vert \mathcal{O}(\beta^2) +R \nonumber \\
&= C^{\ast}-C_{a,0}-\tilde{C}_{-a} + R + \delta_{a,T}(\beta) \nonumber \\
&\leq (C^{\ast}-C_{a,0}-\tilde{C}_{-a}) +\min_{a \in \mathcal{A}}\{C_{a,0}+\tilde{C}_{-a}-C^{\ast}\}+\delta_{a,T}(\beta)
\\& \leq \delta_{a,T}(\beta)
\end{align*}
where we define 
\begin{equation*}
\delta_{a,T}(\beta):=C_{a,T+1}-C^{\ast }+C_{-a,T+1}+\tilde{C}_{-a}-\tilde{C}_{-a,T+1}+  2T\left\vert \mathcal{A}-1\right\vert \mathcal{O}%
(\beta ^{2})
\end{equation*}
First inequality comes from the upper bound on $R$ stated at \eqref{eq_participation_fee_prime}. Second inequality follows from the fact that $(C^{\ast}-C_{a,0}-\tilde{C}_{-a}) +\min_{a \in \mathcal{A}}\{C_{a,0}+\tilde{C}_{-a}-C^{\ast}\}<0$.  
Note that Theorem \ref{thm_convergence} also applies in the case that area $a$ is excluded from the coupling. This implies that $ \tilde C_{-a,T+1} \rightarrow \tilde C_{-a}$, wherein $\tilde C_{-a}$ is the optimal cost to the centralized DC-OPF problem when area $a$ is excluded from the coupling. Hence, by Theorem 1, $\delta_{a,T}(\beta
)\rightarrow 0$ as $\beta \rightarrow 0^{+\text{ }}$and $T\rightarrow
\infty $.
\end{proof}

\medskip 
In Proposition \ref{prop_rational}, we analyze the incentive to participate in coupling using the net cost reduction at the end of the iterative mechanism. While at each iteration some areas may lose money, no area will lose more than $\delta_{a,T}(\beta)$ amount by time $T$ with $\delta_{a,T(\beta)}$ shrinking to zero as $\beta \rightarrow 0^{+\text{ }}$and $T\rightarrow
\infty $.

\subsection{Incentive Compatibility}

\begin{theorem} \label{thm_Nash}
Assume the reported supply costs functions for internal market clearing are $C_{a,g}(\cdot ), g \in \mathcal{G}_a, a \in \mathcal{A}$. Any unilateral deviation reporting of the form $C_{a,g}^{D}(\cdot ) \neq C_{a,g}(\cdot )$ results in a total net cost reduction of at most $\epsilon(\beta,T)$ with $\epsilon(\beta,T) \rightarrow 0$ as $\beta\rightarrow 0^{+}, T \rightarrow \infty$.
\end{theorem}

\begin{proof}
Consider area $a$ reporting information consistent with strictly convex differentiable cost functions $C_{a,g}^{D}(\cdot )$ that differs from the reported supply cost functions for internal market clearing $ C_{a,g}(\cdot )$, i.e. $C_{a,g}^{D}(\cdot ) \neq C_{a,g}(\cdot )$.
Such deviation yields a mechanism output denoted by 
$\boldsymbol{\hat{x}}^{D,k}=(\mathbf{\hat{T}}^{D,k},\boldsymbol{\hat{\theta}}^{D,k},\boldsymbol{\hat{\alpha}}^{D,k})$. Let us denote by $C^F_{a,k}$ the cost for area $a$ after $k$ iterations and $r^{F}_{a,k}$ the intertie flow payments with respect to reported supply costs for internal market clearing, 
\begin{align*}
C_{a,k}^{F}&:=\sum_{g \in \mathcal{G}_{a}}C_{a,g}(\hat{P}_{a,g}^{D,k})  \nonumber\\
r_{a,k}^{F}&:=\hspace{-6pt}\sum_{(i,j)\in \mathcal{T}_{a}}\hspace{-6pt}{\hat{\alpha}_{a^{\prime },j}^{D,k-1}}%
\hat{T}_{a,(i,j)}^{D,k}-\hspace{-6pt}\sum_{(i,j)\in \mathcal{T}_{a}}\hspace{-6pt}\frac{\hat{\mu}
_{(i,j)}^{D,k-1}}{2}(\left\vert \hat{T}_{a,(i,j)}^{D,k}\right\vert -\bar{T}%
_{a,(i,j)}).
\end{align*}
The net cost for area $a$ after $k$ iterations with respect to equilibrium supply costs for internal market clearing is, 
$V^{F}_{a,k}= C^F_{a,k} -r^{F}_{a,k}$. We also define $C^{F}_{-a,k}$ as the total cost of all areas except area $a$ at iteration $k$ when area $a$ deviates from the equilibrium supply function for internal market clearing.

Upon stopping after $T$-iterations the total change in cost inclusive of
incentive transfers and participation fee for area $a$ is%
\begin{align}\label{eq_untruthful_prime}
\sum_{k=0}^{T}\big[V_{a,k+1}^{F}-&V_{a,k}^{F}+\widehat{\Delta \pi} _{a,k}^{F}\big]+R
= (C^{F}-(C_{a,0}+\tilde{C}_{-a})) + R + \delta_{a,T}^{F}(\beta)
\end{align}where 
\begin{equation*}
\delta_{a,T}^{F}(\beta):=C_{T+1}^{F}-C^{F }+\tilde{C}_{-a}-\tilde{C}_{-a,T+1} + 2T\left\vert \mathcal{A}-1\right\vert \mathcal{O}(\beta^{2}).  
\end{equation*}
Similarly, we can write the total change in cost from reporting consistent with internal market clearing as follows%
\begin{align}\label{eq_truthful_prime}
\sum_{k=0}^{T}[V_{a,k+1}&-V_{a,k}+\widehat{\Delta \pi} _{a,k}]+R 
=(C^{\ast }-(C_{a,0}+\tilde{C}_{-a}))+ R+ \delta_{a,T}(\beta) 
\end{align}where
\begin{equation*}
\delta_{a,T}(\beta) :=  C_{T+1}-C^{\ast }+\tilde{C}_{-a}-\tilde{C}_{-a,T+1}     +2T\left\vert \mathcal{A}-1\right\vert \mathcal{O}(\beta ^{2}).
\end{equation*}
Hence by (\ref{eq_untruthful_prime}) and (\ref{eq_truthful_prime}), 
\begin{align*}
\sum_{k=0}^{T}[V_{a,k+1}^{F}-V_{a,k}^{F}+\widehat{\Delta \pi}
_{a,k}^{F}]+R &-(\sum_{k=0}^{T}[V_{a,k+1}-V_{a,k}+\widehat{\Delta \pi}
_{a,k}]+R) \\
&=(C^{F}-C^{\ast })%
+4T\left\vert \mathcal{A}-1\right\vert \mathcal{O}(\beta ^{2}) + (\delta^{F}_{a,T}(\beta)-\delta_{a,T}(\beta)) \\&\geq \epsilon(\beta ,T).
\end{align*}%
The inequality follows from Theorem 1 (optimality of coupling markets in internal equilibrium) and 
\begin{equation*}
\epsilon(\beta ,T)=4T\left\vert \mathcal{A}-1\right\vert \mathcal{O}%
(\beta ^{2}) + \inf_{C_{a,g}^{F}(\cdot )}(\delta^{F}_{a,T}(\beta)-\delta _{a,T}(\beta)).
\end{equation*}Here again, by convergence of the mechanism's output it follows that
$\inf_{C_{a,g}^{F}(\cdot )}(\delta^{F}_{a,T}(\beta)-\delta _{a,T}(\beta)) \rightarrow 0 $ as $T\rightarrow 0$ and $\epsilon(\beta ,T)\rightarrow 0$ as $\beta
\rightarrow 0^{+}$ and $T\rightarrow \infty $.
\end{proof}

\medskip
\noindent
\begin{corollary}
Assume any excess cost reduction for an area (after compensating generators with increased costs) is allocated to consumers.
If the reported supply costs functions for internal market clearing are a Nash equilibrium for each market the original internal equilibrium remains an approximate Nash equilibrium.
\end{corollary}

\begin{proof}
By Theorem 2, the incentive to manipulate the reported supply or demand function is bounded by $\epsilon(\beta,T)$ which is negligible for large $T>0$ and small $\beta>0$.
\end{proof}

\subsection{Budget Balance}

Finally we show the mechanism obtains a surplus that is bounded below by the value of trade (congestion rent) at the tielines under a sufficient condition on the participation fee $R>0$.

%%%%%%%%%%%%%%%%%%%%%%%%%%%%%%%%%%%%%%%%%%%%%%%%%%%%%%%%%%%
%%%%%%%% BEGIN THEOREM %%%%%%%%%%%%%%%%%%%%%%%%%%%%%%%%%%%%%%
%%%%%%%%%%%%%%%%%%%%%%%%%%%%%%%%%%%%%%%%%%%%%%%%%%%%%%%%%%%
\begin{theorem} \label{thm_balance}

We define the mechanism's approximate net budget after $T$ iterations as the sum of all approximate incentive payments \eqref{incentive_hat}, and participation fees
\begin{equation}\label{eq_budget}
\widehat{B}(T,\beta):= \lvert A\rvert {R} +\sum\limits_{a\in \mathcal{A}}\sum\limits_{k=0}^{T}\widehat{\Delta \pi}_{a,k}.
\end{equation}

If $R \geq \min_{a \in \mathcal{A}}\{\tilde{C}_{-a}+C_{a,0}-C^{\ast}\}$
and
\begin{equation}
\label{efficient_coupling_prime}
C_0-C^{\ast} \geq \sum\limits_{a\in \mathcal{A}}[(\tilde{C}_{-a}+C_{a,0}-C^{\ast })-%
\min_{a}\{\tilde{C}_{-a}+C_{a,0}-C^{\ast }\}],
\end{equation} 
the approximate net budget is bounded below by the sum of all trades,
\begin{equation}
\label{eq_budget_balance_hat}
    \lim_{T\to \infty, \beta \to 0^+} \widehat{B}(T,\beta) \geq \sum_{(i,j)\in \mathcal{T}}(\alpha _{a^{\prime },j}^*-\alpha _{a,i}^*) {T}_{a,(i,j)}^*,
\end{equation}
where $\ccalT := \bigcup_{a\in\ccalA} \ccalT_a$ denotes the set of all tielines.
\end{theorem}
%%%%%%%%%%%%%%%%%%%%%%%%%%%%%%%%%%%%%%%%%%%%%%%%%%%%%%%%%%%
%%%%%%%% END THEOREM %%%%%%%%%%%%%%%%%%%%%%%%%%%%%%%%%%%%%%
%%%%%%%%%%%%%%%%%%%%%%%%%%%%%%%%%%%%%%%%%%%%%%%%%%%%%%%%%%%

\begin{proof}
See appendix.
\end{proof}

\begin{remark}
 {\bf (Efficiency vs Information Rents)} Condition \eqref{efficient_coupling_prime} can be interpreted as a requirement on the efficiency of coupling all areas in $\mathcal{A}$.
The right hand side is the total sum of individual surpluses or sum total of information rents (i.e. the cost entailed by coupling all areas without affecting internal equilibrium). The left hand side is the total system benefit from coupling $C_0-C^{\ast}$.
Therefore, if the inequality does not hold then the total information rent (i.e. costs of coupling all areas in $\mathcal{A}$ without affecting internal equlibrium) {\em exceeds} the benefits from coupling.
\end{remark}
\begin{remark}
{\bf (Participation vs Incentive Compatibility):}
Proposition 1 and Theorem 3 characterize the trade-offs between a lower participation fee (which ensures all areas benefit from coupling) versus a higher participation fee (which ensures there is no deficit involved when coupling all areas). Only when the participation fee is set to equal the minimum individual net cost reduction accross all areas, i.e. $R=\min_{a \in \mathcal{A}}\{\tilde{C}_{-a}+C_{a,0}-C^{\ast}\}$ the mechanism can offer non-negative net cost reductions for all areas and incur no deficit. 
\end{remark}

\subsection{A Toy Example}
\begin{figure}[h] 
	\centering
\vspace{-5pt}
	\includegraphics[width=0.5\linewidth]{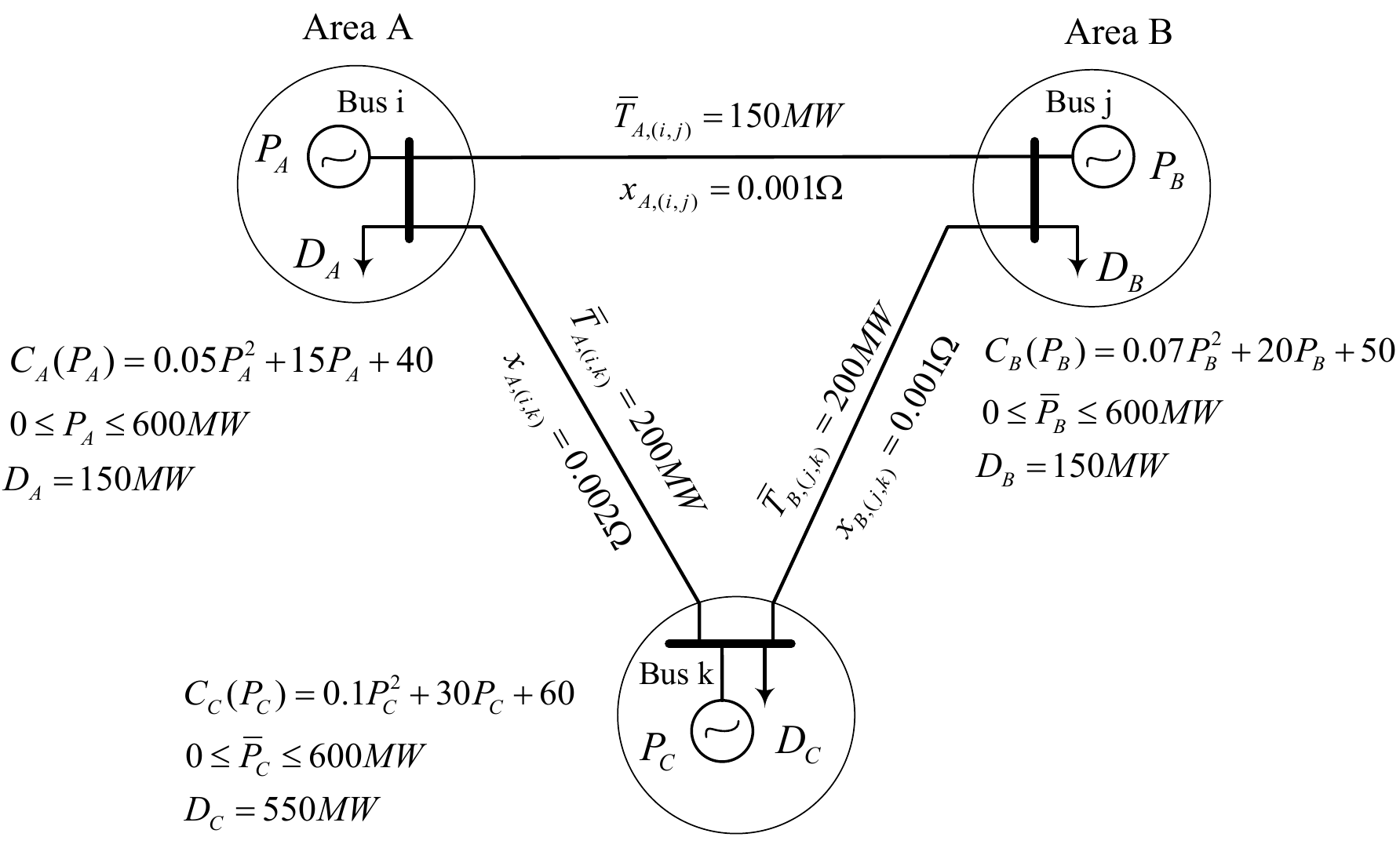}
    \vspace{-5pt}
	\caption{\small Three-bus test system}
       \label{fig.ex}
\end{figure}

Suppose the abstract test system in Fig. \ref{fig.ex}, incorporating three areas A, B, and C, with the generation cost functions, power generation limits, and loads specified below/beside each area, and power flow capacity and reactances specified above and below the tielines.
We show the operation costs of areas in Table \ref{Table.CostEx} for five coupling scenarios, ranging from full coupling to independent operation of areas, and furnish the optimal coupling decisions in Table \ref{Table.TL_LMP_Ex}. 

\begin{table*}[h]
\centering
\caption{Operation costs of areas w.r.t coupling scenarios for the three-bus test system (the symbols ``+"/``," between areas indicate whether the associated areas are coupled/decoupled)}
\begin{tabular}{c c c c c c}
\cline{2-6}
                                           & \multicolumn{5}{c}{Coupling Scenarios}                                                                                                        \\ \cline{2-6} 
                                           & (A+B+C)                    & A,(B+C)                    & (A+C),B                    & (A+B),C                    & A,B,C                      \\ \hline
\multicolumn{1}{c}{Operation Cost of A}  & 15,568                      & 3,415                       & 11,415                      & 4,895                       & 3,415                       \\ %\hline
\multicolumn{1}{c}{Operation Cost of B}  & 8,454                       & 15,625                      & 4,625                       & 2,893                       & 4,625                       \\ %\hline
\multicolumn{1}{c}{Operation Cost of C}  & 9,436                       & 22,810                      & 22,810                      & 46,810                      & 46,810                      \\ \hline
\multicolumn{1}{c}{Total Operation Cost} & 33,457 & 41,850 & 38,850 & 54,598 & 54,850 \\ \hline
\end{tabular}
\label{Table.CostEx}
\end{table*}

\begin{table*}[h]
\centering
\caption{Optimal coupling decisions of the three-bus test system: locational marginal prices, tieline power flows, and capacity prices}
\label{Table.TL_LMP_Ex}
\begin{tabular}{c c c c c c c c c}
\hline
\multicolumn{3}{c}{\begin{tabular}[c]{@{}c@{}}Locational Marginal \\ Prices (\$/MWh)\end{tabular} } & \multicolumn{3}{c}{\begin{tabular}[c]{@{}c@{}}Power Flows \\ (MW)\end{tabular}}         & \multicolumn{3}{c}{\begin{tabular}[c]{@{}c@{}}Capacity Prices \\ (\$/MWh)\end{tabular} } \\ \hline
$\alpha_{A,i}$    & $\alpha_{B,j}$    & $\alpha_{C,k}$    & $T_{A,(i,j)}$ & $T_{A,(i,k)}$ & $T_{B,(j,k)}$ & $\mu_{(i,j)}$ & $\mu_{(i,k)}$ & $\mu_{(j,k)}$ \\ \hline
57.7              & 52.5              & 68.2              & 118.1         & 159.0         & 200.0         & 0             & 0             & 21.0          \\ \hline
\end{tabular}
\end{table*}

\subsubsection{Sensitivity Analysis w.r.t Participation Fee}
Given the data in Table \ref{Table.CostEx}, we report the marginal contributions of areas to coupling and their total cost reductions in Table \ref{Table.Ex1}, for the exact estimation of the participation fee (i.e., $R=\min_{a \in \mathcal{A}}\{\tilde{C}_{-a}+C_{a,0}-C^{\ast}\}$), as well as its overestimation and underestimation.  
The exact value of the participation fee equals to $\$5,393$, as shown in the first row of Table \ref{Table.Ex1}, while the total cost reduction of each area amounts to the sum of participation fee, marginal contribution to coupling, and the internal cost saving of that area.
Since none of the areas incur any loss, the exact participation fee corroborates the individual rationality of the proposed mechanism.
The congestion rent in Table \ref{Table.CostEx} represents the net monetary value of all tieline transactions, given by the term $\sum_{(i,j)\in \mathcal{T}}(\alpha _{a^{\prime },j}^*-\alpha _{a,i}^*) {T}_{a,(i,j)}^*$ in the right hand side of \eqref{eq_budget_balance_hat}, which we calculate using the data in Table \ref{Table.TL_LMP_Ex}.
In addition, the mechanism's surplus indicates the net budget of the mechanism (defined as $B(T)$ in \eqref{eq_budget}) minus the congestion rent. 
Note that we subtract the congestion rent from the net budget as it has a different claimant.  %, while securing a positive surplus.

Now, imagine the market maker overestimates or underestimates the participation fee by $\$2,500$, as we show in the second and third rows of Table \ref{Table.Ex1}, respectively.
In the former case, area B incurs a loss of $\$2,500$ as a result of participating in market coupling, and probably opts out, while in the latter case all areas benefit from participating in the market coupling, yet the mechanism leads to a deficit of $\$4,855$.
It is worth highlighting that, in the cases of exact estimation and overestimation of participation fees both conditions of the Theorem \ref{thm_balance} hold, thus, the mechanism's surplus ($B(T)$ minus congestion rent) is non-negative as anticipated by Theorem \ref{thm_balance}. 
However, in the case of underestimated participation fee, the first condition of Theorem \ref{thm_balance} is violated leading to a negative surplus and violation of budget balance. 

\begin{table*}[h]
\centering
\caption{Summary of areas cost reduction w.r.t changes in participation fee}
\begin{tabular}{c c c c c c c c c c}
\cline{3-10}
                                              &       & \multicolumn{3}{c}{\begin{tabular}[c]{@{}c@{}}Marginal Contributions \\ to Coupling (\$)\end{tabular}}                        & \multicolumn{3}{c}{\begin{tabular}[c]{@{}c@{}}Total Cost  \\ Reduction (\$)\end{tabular}} & \begin{tabular}[c]{@{}c@{}}Congestion \\ Rent (\$)\end{tabular}
                 &  \begin{tabular}[c]{@{}c@{}}Mechanism's \\Surplus  (\$)\end{tabular}
                 \\ \cline{3-8}
                                              &       & A                        & B                       & C                       & A               & B              & C              &  &                                        \\ \hline
\multicolumn{1}{c}{} & 5,393 & -20,545                     &           -9,222              &       16,234 &     3,000      &    0      &   15,748       &  4,193 & 2,645                                \\ %\cline{2-2} \cline{6-9} 
\multicolumn{1}{c}{\begin{tabular}[c]{@{}c@{}}Participation \\ Fee (\$)\end{tabular}}                        & 7,893 &     -20,545                     &           -9,222              &       16,234                  & 500                & -2,500          &    13,248      &     4,193 & 10,145                             \\ 
%\cline{2-2} \cline{6-9} 
%\multicolumn{1}{c}{}   
& 2,893 & -20,545                     &           -9,222              &       16,234       & 5,500          & 2,500         & 18,248        &  4,193 &  -4,855                                \\ \hline
\end{tabular}
\label{Table.Ex1}
\end{table*}

%There is a critical participation fee, equal to $R=5,393-\frac{6,838}{3}=\$3,114$, for which the mechanism leads to a zero net budget and all areas benefit from participating in the market coupling, while any participation fee below this critical value would lead to violating the budget balance and incurring a deficit. 

\subsubsection{The Strategic Effects of Coupling}
Here we aim to investigate how a deviation from equilibrium in area A impacts its total cost reduction.
To this end we consider two cases: Case 1 where area A uses a cost function greater than its equilibrium cost function by a factor of $1.1$, i.e., $C^{D}_A(P_A)=1.1C_A(P_A)$; Case 2 where area A uses a cost function lower than its equilibrium cost function by a factor of $0.9$, i.e., $C^{D}_A(P_A)=0.9C_A(P_A)$.
We present the cost reduction analysis of the two cases in Table \ref{Table.Ex2} for the participation fee equal to $\$5,393$.
As shown in Table \ref{Table.Ex2}, and compared to the first row of Table \ref{Table.Ex1}, any deviation from equilibrium reporting leads to a lower cost reduction for Area A, thus, area A is better off with equilibrium reporting. 
In addition, deviation from equilibrium reporting of area A may cause other areas undergo loss and opt out, such as area B in Case 1, or may lead to reducing the mechanism's surplus compared to equilibrium reporting, as in Case 2.

\begin{table*}[h]
\centering
\caption{Summary of areas' net cost reduction for deviation in area A. \newline Cases 1 and 2 respectively refer to overbidding and underbidding}
\begin{tabular}{c c c c c c c c c c}
\cline{2-10}
                             & \begin{tabular}[c]{@{}c@{}}Participation \\ Fee (\$)\end{tabular} & \multicolumn{3}{c}{\begin{tabular}[c]{@{}c@{}}Marginal Contributions \\ to Coupling (\$)\end{tabular}} & \multicolumn{3}{c}{\begin{tabular}[c]{@{}c@{}}Total Cost \\ Reduction (\$)\end{tabular}} &  \begin{tabular}[c]{@{}c@{}}Congestion \\ Rent (\$)\end{tabular}             
                 &            \begin{tabular}[c]{@{}c@{}}Mechanism's  \\Surplus  (\$)\end{tabular} \\ \cline{3-8}
                             &                                         & A                   & B                 & C                  & A          & B         & C          &                                         \\ \hline
\multicolumn{1}{c}{Case 1} & 5,393                                   & -18,795             & -9,958            & 16,555            & \textbf{2,919}     & -330       & 14,742    & 3,982 &3,981                                   \\ %\hline
\multicolumn{1}{c}{Case 2} & 5,393                                   & -22,408             & -8,609            & 15,837            & \textbf{2,899}     & 511     & 16,884    & 4,428 & 998                                     \\ \hline
\end{tabular}
\label{Table.Ex2}
\end{table*}

\vspace{-5pt}

\section{Numerical Testbed} \label{sec.results}
Here we provide the numerical results for the three-zone IEEE Reliability Test System (IEEE-RTS) with 96 generating units, 73 buses,  115 transmission lines, and 5 tielines.
The topology of the test system under study, the operation limitations of generating units and lines/tielines, as well as the nodal load data are available in \cite{grigg1999ieee}.
We make the following changes to the original network in order to simulate the occasion where transmission lines/tielines are congested: 
\begin{itemize}
\item The capacity of transmission line located in area 1 which connects the buses $(16,17)$ is reduced from 500MW to 200MW, e.g., $\overline{F}_{1,(16,17)}$=200MW.
\item The capacity of transmission line located in area 2 which connects the buses $(3,24)$ is reduced from 400MW to 150MW, e.g., $\overline{F}_{2,(3,24)}$=150MW.
\item The capacity of tieline 4 is reduced from 500MW to 100MW.
\end{itemize}
For quick accessibility, we also provide the tieline data in Table \ref{table.TLData}, where the set of areas is defined as $\mathcal{A}=\{\text{A},\text{B},\text{C}\}$, in accordance with \cite{grigg1999ieee}.
%\begin{table} [ht]
% \vspace{-5pt}
%	\centering
%	\caption{\footnotesize Tieline Data}
%    \vspace{5pt}
%	\includegraphics[width=0.9\linewidth]{TLData.pdf}
%    \vspace{-5pt}
%	\label{table.TLData} 
%\end{table}

\begin{table}[h]
	\centering
	\caption{ Tieline data for the three-area IEEE-RTS}
	\label{table.TLData} 
\begin{tabular}{c c c c c}
\hline
\begin{tabular}[c]{@{}c@{}}Tieline \\ Number\end{tabular} & \begin{tabular}[c]{@{}c@{}}Sending End\\ (Area-Bus)\end{tabular} & \begin{tabular}[c]{@{}c@{}}Receiving End \\ (Area-Bus)\end{tabular} & \begin{tabular}[c]{@{}c@{}}Reactance \\ (pu)\end{tabular} & \begin{tabular}[c]{@{}c@{}}Capacity \\ (MW)\end{tabular} \\ \hline
TL1& A-7& B-3  & 0.16  & 175    \\ 
TL2& A-13 & B-15& 0.08 & 500    \\ 
TL3& A-23& B-17& 0.07& 500      \\ 
TL4& C-25& A-21& 0.10& 100      \\ 
TL5& C-18& B-23& 0.10& 500      \\ \hline
\end{tabular}
\end{table}

We carry out the simulation for a single hour corresponding to the peak load of the three-zone IEEE-RTS, and assume all generating units are  online. 
In addition to the simulation results of the proposed incentive-compatible approach, we also furnish the results associated with \textit{locational marginal price remuneration scheme}, where the areas are merely paid/charged for the exported/imported energy over tielines at the locational marginal prices of boundary buses. 
This remuneration scheme is in accordance with the definition of optimal value in \eqref{eq_optimal_k}, and serves as a benchmark case. 
The cost functions are quadratic.

\subsection{Numerical Results of Iterative Capacity Allocation}
Here we present the numerical results of implementing the proposed iterative capacity allocation method on the three-area IEEE-RTS.
The constant step-size associated with the capacity price updates is $\beta=0.3$, the diminishing step-size associated with intertie power flow updates as well as the voltage phase angles and locational marginal price updates are deemed as $\rho_k=\frac{1}{1+\text{log}(k)}$. 
%while the iterative process terminates as the following conditions are met:
%
%\begin{align*}
%&\left\vert\hat{T}^k_{a,(i,j)}+\hat{T}^k_{a',(i,j)}\right\vert\leq 0.5\\ 
%&\left\vert\mu^k{(i,j)}-\mu^{k-1}{(i,j)}\right\vert\leq 0.1.
%\end{align*}

The simulation results indicate that the values obtained from the iterative capacity allocation converge to that of the centralized DC-OPF model, confirming the efficiency of the proposed method. 

\begin{figure}[h] 
	\centering
\vspace{-5pt}
	\includegraphics[width=0.5\linewidth]{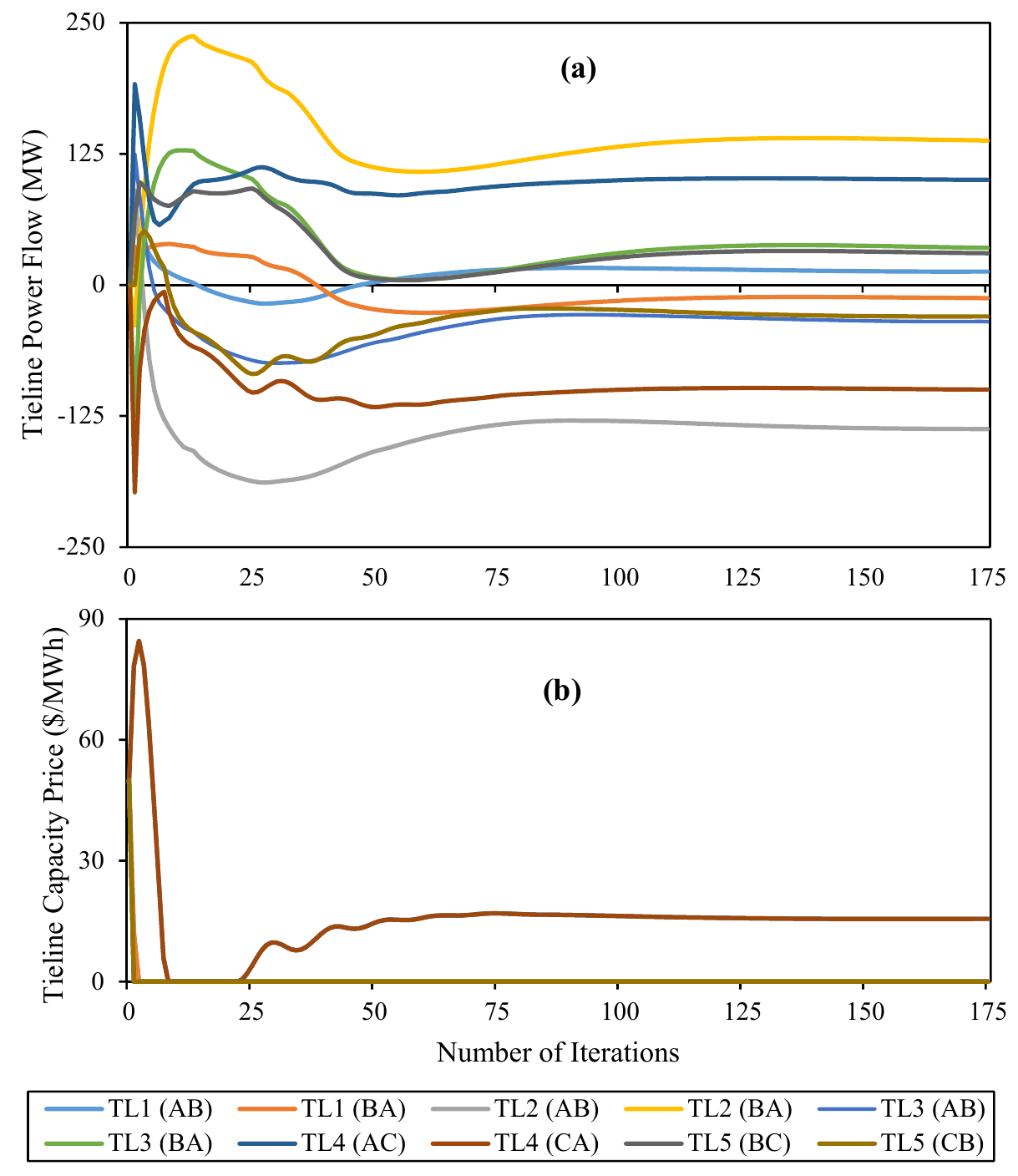}
    \vspace{-10pt}
	\caption{\small (a) Tieline power flows, $T_{a,(i,j)}^{k},$ (MW) (b) Tieline capacity prices, $\mu _{a,(i,j)}^{k},$ (\$/MW) }
       \label{fig.TC}
\end{figure}

\begin{figure}[h] 
	\centering
\vspace{-5pt}
	\includegraphics[width=0.5\linewidth]{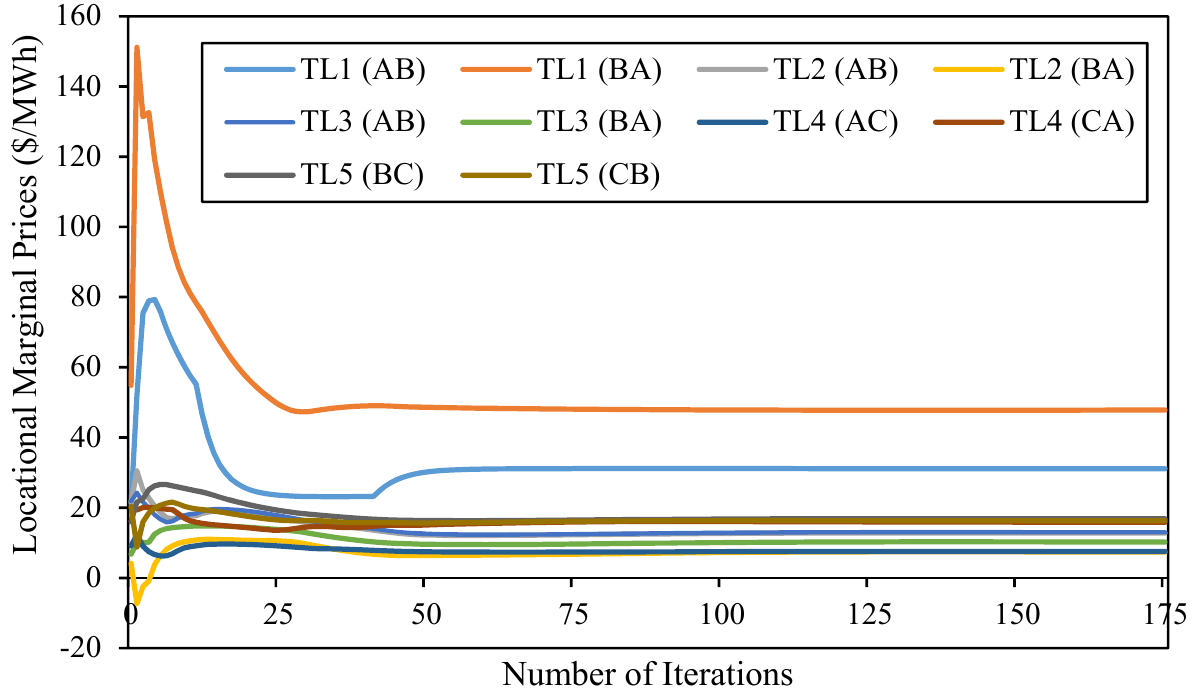}
    \vspace{-5pt}
	\caption{\small Locational marginal prices of tieline incident buses (\$/MWh) }
       \label{fig.ImportExport}
\end{figure}

Figures \ref{fig.TC}-(a) and (b) respectively represent the intertie power flows and the associated capacity prices for each iteration, converging after 175 iterations.
The initial values of intertie power flows are all zero, meaning that the the initial state corresponds to the independent market clearing of all areas.
In addition, the capacity prices in Fig. \ref{fig.TC}-(b) are all initialized at $\$50$ per MWh, which is slightly greater than the highest incremental cost rate of the most expensive generator.
The power flow of tieline number 4 converges to its maximum limit, 100MW, and bears a capacity price of $\$15.6$ per MWh, while the power flows of other tielines remain within their operating limits and all the associated capacity prices converge to zero.
We present the locational marginal prices at tieline incident buses, or indeed the import/export quotes of the areas, in Fig. \ref{fig.ImportExport} for all iterations.
Further, we provide the optimal tieline power flows, capacity prices, and the locational marginal prices at incident buses in Table II, as they are essential to calculating the congestion rent. The congestion rent in Table IV represents the net monetary value of all tieline transactions, given by the term $\sum_{(i,j)\in \mathcal{T}}(\alpha _{a^{\prime },j}^*-\alpha _{a,i}^*) {T}_{a,(i,j)}^*$.
In addition, the mechanism's surplus indicates the net budget of the mechanism (defined as $\widehat{B}(T)$ minus the congestion rent. 
Note that we subtract the congestion rent from the net budget as it has a different claimant. 

%The import/export bids of areas for the energy carried by incident tielines are presented in Fig. \ref{fig.ImportExport}.

\begin{table*}[h]
	\centering
	\caption{ Optimal coupling decisions for the three-area IEEE-RTS: Tieline power flows, capacity prices, and locational marginal prices at incident buses}
	\label{table.TL_LMP_IEEE} 
\begin{tabular}{c c c c c}
\hline
\begin{tabular}[c]{@{}c@{}}Tieline \\ Number\end{tabular}
&
\begin{tabular}[c]{@{}c@{}} Power Flow\\ (MW)\end{tabular}&\begin{tabular}[c]{@{}c@{}}Capacity Price\\ (\$/MWh)\end{tabular}   & \begin{tabular}[c]{@{}c@{}}Sending End LMP\\ (\$/MWh)\end{tabular}   & \begin{tabular}[c]{@{}c@{}}Receiving End LMP\\ (\$/MWh)\end{tabular}     \\ \hline
TL1 & 13.1 &  0           & 31.1                     & 47.8                                            \\ 
TL2     & -136.5  & 0       & 12.7                     & 7.3                                         \\ 
TL3     & -34.3      &  0       & 13.0                     & 10.2                                     \\ 
TL4     & -100.0         & 15.6       & 15.9                     & 7.5                                  \\ 
TL5     & -29.3          &  0       & 16.3                     & 16.8                                 \\ \hline
\end{tabular}
\end{table*}

\subsection{Incentive Transfers Analysis}
Here we first provide the results of truthful reporting and next consider the untruthful behaviour of areas.

\subsubsection{Equilibrium Reporting}
We provide the operation costs of areas in Table \ref{Table.CostIEEE} for the full coupling and independent operation of areas, as well as coupling with one area left out.
%We use the same notation for describing these coupling scenarios as that in the toy example, e.g., the symbol ``+"/``," between areas indicates coupling/decoupling of areas.
As expected, any coupling leads to a reduction in total operation cost of the system compared to independent operation of areas, and the highest system saving, $\$3,866$, is realized for the case of full coupling.

%Also, area A demonstrate the greatest contribution to coupling in terms of cost reduction, while the areas B and C come next.

\begin{table*}[h]
\centering
\caption{Operation costs of areas w.r.t coupling scenarios for the three-area IEEE-RTS (the symbols ``+"/``," between areas indicate whether the associated areas are coupled/decoupled)}
\begin{tabular}{c c c c c c }
\cline{2-6}
                                           & \multicolumn{5}{c}{Coupling Scenarios}                                                                                                        \\ \cline{2-6} 
                                           & (A+B+C)                    & A,(B+C)                    & (A+C),B                    & (A+B),C                    & A,B,C                      \\ \hline
\multicolumn{1}{c}{Operation Cost of A}  & 29,822	& 31,843	&32,757	&28,986	&31,843
                      \\ 
\multicolumn{1}{c }{Operation Cost of B}  & 34,062	&35,295	&33,366	&33,667	&33,366
              \\ 
\multicolumn{1}{c}{Operation Cost of C}  & 25,995	&26,162	&26,534	&28,537	&28,537
              \\\hline 
\multicolumn{1}{c }{Total Operation Cost} & 89,879	&93,300	&92,657	&91,190	&93,745
 \\ \hline
\end{tabular}
\label{Table.CostIEEE}
\end{table*}

\begin{table*}[h]
\centering
\caption{Summary of areas' net cost reduction for equilibrium reporting}
\label{table.CostRedIEEE}
\begin{tabular}{c c c c c c c c c}
\hline
\begin{tabular}[c]{@{}c@{}}Participation \\ Fee (\$)\end{tabular} & \multicolumn{3}{c}{\begin{tabular}[c]{@{}c@{}}Marginal Contributions \\ to Coupling (\$)\end{tabular}} & \multicolumn{3}{c}{\begin{tabular}[c]{@{}c@{}}Total Cost \\ Reduction (\$)\end{tabular}} &
\begin{tabular}[c]{@{}c@{}}Congestion \\ Rent (\$)\end{tabular}
&\begin{tabular}[c]{@{}c@{}}Mechanism's Net\\ Budget (\$)\end{tabular} \\ \cline{2-7}
                                        & A                   & B                  & C                 & A                & B                & C           &                                         \\ \hline
1311                                    & -1,400              & -3,474             & 1,230             & 2110            & 1466            & 0           &  1,872    & 290                                  \\ \hline
\end{tabular}
\end{table*}

Given the information in Table \ref{Table.CostIEEE}, the best estimate of the participation fee equals to $\$1,311$, as shown in Table \ref{table.CostRedIEEE}.
Area B receives the highest reward due to marginal contribution to coupling and the area A comes next; however, the marginal contribution of area C to coupling is a positive number meaning that it is charged rather than being rewarded (Table \ref{table.CostRedIEEE}).
The sum of participation fee, marginal contribution to coupling, and the internal cost saving of each area amounts to its total cost reduction.
The proposed incentive payment approach results in non-negative cost reductions for all areas, where areas A and C respectively benefit the most and least.
Since the participation fee meets the first condition of Theorem \ref{thm_balance}, and the second condition also holds true (sum of area cost reductions is less than the system cost saving, i.e., $\$3,576<\$3,866$), the mechanism bear a non-negative surplus as suggested by Theorem \ref{thm_balance}.

\subsubsection{Deviation Reporting}
In this part we consider three deviation reporting cases. In each case deviations only come from one area at a time.
The deviation cost functions are greater than given cost functions by a factor of $1.1$, meaning that, the areas A, B, and C respectively use the cost functions $C^{D}_{A,g}(P_{A,g})=1.1C_{A,g}(P_{A,g}), g\in\mathcal{G}_A$, $C^{D}_{B,g}(P_{B,g})=1.1C_{B,g}(P_{B,g}), g\in\mathcal{G}_B$, and $C^{D}_{C,g}(P_{C,g})=1.1C_{C,g}(P_{C,g}), g\in\mathcal{G}_C$.
We have provided the total cost reductions of areas in Tables \ref{table.CostRedProp} and \ref{table.CostRedLMP} respectively for the proposed incentive-compatible approach and the locational marginal price remuneration scheme.
These tables contain the results for both deviation and from equilibrium reporting cases.

\begin{table}[h!]
\centering
\caption{Cost reductions of areas for equilibrium reporting and deviation reporting under proposed incentive payment approach}
\label{table.CostRedProp}
\begin{tabular}{c c c c}
\cline{2-4}
& \multicolumn{3}{c}{Total Cost Reduction (\$)} \\ \cline{2-4} & A& B & C\\ 
\hline
\multicolumn{1}{c}{Equilibrium Reporting}      &  $\bf 2,110$      & $\bf 1,466$       &  $\bf 0$
\\ 
\hline
\multicolumn{1}{c}{A Deviates} &  $\bf 2,074$      & $1,796$      &    $101$             \\ 
%\hline
\multicolumn{1}{c}{B Deviates}  &    $2,272$    & $\bf 1,458$      &      $47$    
\\ 
%\hline
\multicolumn{1}{c}{C Deviates} &    $2,089$    &      $1,525$     & $\bf -5$
\\ 
\hline
\end{tabular}
\end{table}

\begin{table}[h!]
\centering
\caption{Cost reductions under deviation and equilibrium reporting under locational marginal price remuneration scheme}
\label{table.CostRedLMP}
\label{table.incentive4}
\begin{tabular}{c c c c}
\cline{2-4}
& \multicolumn{3}{c}{Total Cost Reduction (\$)} \\ \cline{2-4} & A& B & C\\ \hline
\multicolumn{1}{c}{Equilibrium Reporting}      & $\bf 2,886$       &  $\bf 1,556$      & $\bf 1,296$
\\ 
\hline
\multicolumn{1}{c}{A deviates} & $\bf 2,711$       &   $1,669$     & $1,206$             \\ 
%\hline
\multicolumn{1}{c}{B deviates} &     $2,883$  &    $\bf 1,712$    &  $1,223$    
\\ 
%\hline
\multicolumn{1}{c}{C deviates} &   $3,015$        &    $1,571$    &  $\bf 1,298$
\\ 
\hline
\end{tabular}
\end{table}

None of the areas manage to improve their net cost reduction through deviation reporting, and all the attempts lead to lower cost reductions compared to the equilibrium reporting counterpart (Table \ref{table.CostRedProp}).
Since the agent attempting  a manipulation receives a fraction of the net cost reduction resulting from coupling, there is no incentive to manipulate the supply function in the original internal market equilibrium.
%Since areas B and C acquire a greater cost reduction through reporting untruthfully under the locational marginal price remuneration scheme (Table \ref{table.CostRedLMP}). The locational marginal price remuneration scheme is prone to untruthful behaviour of areas, whereas the proposed approach is not.
\begin{remark} {\bf Coupling with LMP-based remuneration}
A coupling scheme based upon locational marginal prices (LMPs) induces incentive compatibility problems. To illustrate this fact, in (Table \ref{table.CostRedLMP}) we describe the resulting net benefits for each area. Individual agents in markets B and C may obtain a rent by manipulating the reported information.
\end{remark}
\section{Conclusions}

We have introduced and analyzed a market design for scheduling and pricing power flows between interconnected electricity markets. Each area operator participates (on behalf of agents active in its market) by iteratively submitting bids for trading energy across interties and the prices for interconnection capacity are adjusted as a function of excess demand. 
The proposed design allows individual area operators to retain local control and can be easily implemented {\em after} existing short-run markets (e.g. hour ahead) have cleared. The market coupling design does not alter the structure of incentives in each internal market, i.e. any internal market equilibrium will remain so (approximately) after coupling is implemented. This is achieved with incentive transfers (updated at each iteration) that remunerate each area with its marginal contribution (i.e. cost savings) to all other participating areas.
Finally, we identify a sufficient condition on a uniform participation fee ensuring the mechanism incurs no deficit.

\section{Appendix}

\subsection{Preliminary results}

\begin{proposition}\label{prop_assumption}
If there exists a feasible solution to the decentralized DC-OPF problem \eqref{eq.Obj2}-\eqref{eq.slack2} for area $a\in \mathcal A$ when $T_{a,(i,j)}=0$ for all $(i,j) \in \mathcal{T}_a,$ then for each intertie $
(i,j)$, there exists a critical price $\bar{\mu}_{(i,j)}>0$ such that for $\mu
_{(i,j)}\geq\bar{\mu}_{(i,j)}$ the optimal solution to \eqref{eq.Obj2} -  \eqref{eq.slack2} will have no flow along the intertie, i.e. $T_{a,(i,j)}=0$. 
\end{proposition}
\begin{proof}
We prove by contradiction. Suppose that there exists an optimal solution such that for a tieline $(i^{*},j^{*}) \in \mathcal{T}_{a}$, we have $T^\epsilon_{a,(i^{*},j^{*})}=\epsilon \neq 0$ for any capacity price ${\mu}_{(i^*,j^*)}>0$. Denote the power generation and tieline flows of this optimal solution with $P_{a,g}^\epsilon$ and $T^\epsilon_{a,(i,j)}$, respectively. In this case, the optimal objective value in \eqref{eq.Obj2} is given by 
\begin{align}\label{eq:max-intertie}
\sum_{g\in \mathcal{G}_a} C_{a,g}(P^{\epsilon}_{a,g}) &- \epsilon\alpha_{a',j^{*}} 
+\frac{{\mu}_{(i^{*},j^{*})}}{2}(|\epsilon|-\bar{T}_{a,(i*,j*)})   - \sum_{(i,j)\in \mathcal{T}_{a} \setminus (i^{*},j^{*}) }\frac{\mu_{(i,j)}}{2}\bar{T}_{a,(i,j)}. 
\end{align}
Consider a feasible solution with $T_{a,(i,j)}=0$ for all $(i,j) \in \mathcal{T}_a$. We denote the power generation of the feasible solution with zero tieline flows by $P_{a,g}^0$. Given the feasible solution, the objective value is given by $\sum_{g\in \mathcal{G}_a} C_{a,g}({P}^0_{a,g}) -\sum_{(i,j)\in \mathcal{T}_a}\frac{{\mu}_{(i,j)}}{2} \bar{T}_{a,(i,j)}.$ \par
Considering the difference between the objective value in \eqref{eq:max-intertie} with that of the feasible solution when $T_{a,(i^*,j^*)}=0$, we have
\begin{equation}\label{eq:max-intertie-result}
\sum_{g\in \mathcal{G}_a} C_{a,g}(P^{\epsilon}_{a,g})-\epsilon \alpha_{a',j^{*}} +\frac{{\mu}_{(i^{*},j^{*})}}{2}|\epsilon|  - 
\sum_{g\in \mathcal{G}_a} C_{a,g}({P}^0_{a,g}) <0
\end{equation}
which must be negative by the assumption that  $(P_{a,g}^\epsilon,T_{a,(i,j)}^\epsilon)$ is optimal. Observe now that there exists a capacity price value $\bar{\mu}_{(i^{*},j^{*})}>0$ such that the difference in objective values \eqref{eq:max-intertie-result}  becomes positive for any $\epsilon\neq 0$. Thus, for each intertie $%
(i,j)$, there exists a price $\bar{\mu}_{(i,j)}>0$ such that for $\mu
_{(i,j)}\geq \bar{\mu}_{(i,j)}$ the optimal solution will have no flow along the intertie, i.e. $T_{a,(i,j)}=0$.
% \red{\{RK: It seems to me that both $\epsilon>0$ and $\epsilon<0$ should be covered, meaning that for a very large capacity price the area $a$ is neither willing to sell nor willing to buy.
% Maybe we can discuss the latter case in terms of the willingness of the neighbouring area to sell, but we at least need to discuss it.\}}
%
%
% While we are increasing $\overline{\mu}_{(i^{*},j^{*})},$ the inequality in (\ref{eq:max-intertie-result}) does not hold after some level. This contradicts with the assumption of there exists a feasible solution with same capacity price $\overline{\mu}_{(i^{*},j^{*})}$ such that there exists $(i^{*},j^{*}) \in \mathcal{T}_{a}$ with $T_{a,(i^{*},j^{*})}=\epsilon > 0.$ Therefore, for each intertie $%
% (i,j)$, there exists a price $\bar{\mu}_{(i,j)}>0$ such that for $\mu
% _{(i,j)}\geq \bar{\mu}_{(i,j)}$ the optimal solution to \eqref{eq.Obj2} -  \eqref{eq.slack2} will have no flow along the intertie, i.e. $T_{a,(i,j)}=0$. 
\end{proof}

\subsection{Proof of Lemma 2}
Suppose the sequence $\{\mu _{(i,j)}^{k}:k>0\}$ does not converge. So either {\em (i)},
$\mu _{(i,j)}^{k}\rightarrow +\infty $ or {\em (ii)}, it oscillates with $\liminf \mu
_{(i,j)}^{k}<\limsup \mu _{(i,j)}^{k}$. Case {\em (i)} is easily discarded since from the assumption on maximum intertie capacity price:
\[
\mu _{(i,j)}^{k}\rightarrow +\infty 
~~~\Rightarrow~~~
\left\vert T_{a,(i,j)}^{k}\right\vert
+\left\vert T_{a^{\prime },(i,j)}^{k}\right\vert \rightarrow 0
\]
which is a
contradiction.
Now let us consider case {\em (ii)}. Capacity price oscillation implies that for some $\epsilon >0$ 
\begin{align*}
&\lim \sup_{k}\{\frac{\left\vert T_{a,(i,j)}^{k}\right\vert +\left\vert
T_{a^{\prime },(i,j)}^{k}\right\vert }{2}-\bar{T}_{a,(i,j)}\}\geq \epsilon
>0 \geq \lim \inf_{k}\{\frac{\left\vert T_{a,(i,j)}^{k}\right\vert
+\left\vert T_{a^{\prime },(i,j)}^{k}\right\vert }{2}-\bar{T}_{a,(i,j)}\}
\end{align*}

Thus the sequence $\{\frac{\left\vert T_{a,(i,j)}^{k}\right\vert +\left\vert
T_{a^{\prime },(i,j)}^{k}\right\vert }{2}-\bar{T}_{a,(i,j)}\}$ exhibits
infinitely many upcrossings of $\frac{\epsilon }{4}$ and $\frac{3\epsilon }{4%
}$. Let $k>0$ denote an index for an iteration in which an upcrossing of $%
\frac{\epsilon }{4}$ takes place, i.e.:% 
\begin{align}
\frac{\left\vert T_{a,(i,j)}^{k}\right\vert +\left\vert T_{a^{\prime
},(i,j)}^{k}\right\vert }{2}-\bar{T}_{a,(i,j)}&\geq \frac{\epsilon }{4} , \textrm{ and } 
\frac{\left\vert T_{a,(i,j)}^{k-1}\right\vert +\left\vert T_{a^{\prime
},(i,j)}^{k-1}\right\vert }{2}-\bar{T}_{a,(i,j)}<\frac{\epsilon }{4}.%
\end{align}
Since $T_{a,(i,j)}^{k+1}-T_{a,(i,j)}^{k}=\rho _{k}(\hat{T}_{a}^{k}-T_{a}^{k})
$ it follows that 
\begin{align}
\sum_{\ell =1}^{\tau (k)}(T_{a,(i,j)}^{k+\ell }-T_{a,(i,j)}^{k+\ell
-1})&=\sum_{\ell =1}^{\tau (k)}\rho _{k+\ell }(\hat{T}_{a}^{k+\ell
}-T_{a}^{k+\ell })\leq 2\bar{T}\sum_{\ell =1}^{\tau (k)}\rho _{k+\ell }
\end{align}
where $\bar{T}<\infty$ us a uniform upper bound on unconstrained intertie flows. This allows us to obtain a lower bound on number of iterations $\tau (k)>0$
needed to upcross $\frac{3\epsilon }{4}$ from $\frac{\epsilon }{4}$: 
\begin{equation*}
2\bar{T}\sum_{\ell =1}^{\tau (k)}\rho _{k+\ell }>\frac{3\epsilon }{4}-\frac{%
\epsilon }{4}=\frac{\epsilon }{2}
\end{equation*}%
Since %
 $\tau (k)\rho _{k}>\sum_{\ell =1}^{\tau (k)}\rho _{k+\ell }
$,
it follows that $\tau (k)>\frac{2\bar{T}}{\rho _{k}}\rightarrow +\infty $
since $\rho _{k}\rightarrow 0^{+}$. We also have %
$
\mu _{(i,j)}^{k+\tau }\geq \mu _{(i,j)}^{k}+ \beta\tau\frac{\epsilon }{4}
$ %
for all $0<\tau \leq \tau (k)$. Since $
\sum_{\ell =1}^{\tau (k)}\rho _{k+\ell }\rightarrow 0^{+}
$, 
it follows that there exists $\tau >0$ such that $\mu _{(i,j)}^{k+\ell }>%
\bar{\mu}_{(i,j)}$ and $\hat{T}_{a}^{k+\ell }=0$ for all $\tau <\ell \leq
\tau (k)$. Hence, $\left\vert T_{a,(i,j)}^{k+\ell }\right\vert \leq
\left\vert T_{a,(i,j)}^{k+\tau }\right\vert $for all $\tau <\ell \leq \tau
(k)$ which is a contradiction to upcrossing $\frac{3\epsilon }{4}$. 

The non-negativity of $\mu_{(i,j)}^\infty$ follows by the max operation in \eqref{eq.PriceUpdate}. 

\subsection{Proof of Theorem \ref{thm_balance}.}

By substituting (\ref{incentive_hat}) for $%
\widehat{\Delta \pi} _{a,k}$, the mechanism's approximate net revenue comprised can be written as : 
\begin{equation}
\widehat{B}(T,\beta)= \left\vert \mathcal{A}\right\vert R+\sum\limits_{a\in \mathcal{A}}\sum_{k=0}^{T} \widehat{M_{a,k}}+r_{a,k+1}-r_{a,k}.
\end{equation}%
Telescoping sum over $k$ for the second term on the right hand side yields
\begin{align}
\widehat{B}(T,\beta) &= \;\left\vert \mathcal{A}\right\vert R  + \sum_{a\in\ccalA}[ C_{-a,T+1}-C_{-a,0} -(\tilde C_{-a,T+1}-\tilde C_{-a,0})\nonumber \\& \hspace{6pt} + r_{a,T+1}-r_{a,0}]+ 2T\left\vert \mathcal{A}-1\right\vert \mathcal{O}(\beta ^{2})\\
&= \left\vert \mathcal{A}\right\vert R+\sum\limits_{a\in \mathcal{A}}[C_{T+1}-%
\tilde{C}_{-a,T+1}-C_{a,0}+r_{a,T+1}-r_{a,0}]\nonumber \\& \hspace{6pt} +\sum\limits_{a\in \mathcal{A}}(C_{a,0}-C_{a,T+1}) +2T\left\vert \mathcal{A}-1\right\vert \mathcal{O}(\beta ^{2})
\end{align}
where the second equality follows by adding and subtracting  $\sum_{a\in\ccalA} C_{a,T+1}$ and $\sum_{a\in\ccalA} C_{a,0}$, and by using $C_{T+1} = C_{a,T+1} + C_{-a,T+1}$ and $C_{a,0} = \tilde C_{a,0}$. We get the following lower bound by using the minimum value of the participation fee
\begin{align}
  \widehat{B}(T,\beta) &\geq \left\vert \mathcal{A}\right\vert \min_{a}\{\tilde{C}_{-a}+C_{a,0}-C^{\ast}\} \nonumber\\&+\sum\limits_{a\in \mathcal{A}}[C_{T+1}-%
\tilde{C}_{-a,T+1}-C_{-a,0}+r_{a,T+1}-r_{a,0}] 
\nonumber\\&+\sum\limits_{a\in \mathcal{A}}(C_{a,0}-C_{a,T+1}) + 2T\left\vert \mathcal{A}-1\right\vert \mathcal{O}(\beta ^{2}).
\end{align}
We reorganize the terms to get 
\begin{align}
   \widehat{B}(T,\beta)&\geq \sum\limits_{a\in \mathcal{A}}[C_{T+1}-(\tilde{C}_{-a,T+1}+C_{a,0})\nonumber \\& \hspace{6pt}+%
\min_{a}\{\tilde{C}_{-a}+C_{a,0}-C^{\ast}\}]+\sum\limits_{a\in \mathcal{A}%
}(r_{a,T+1}-r_{a,0})
\nonumber \\&\hspace{6pt}+C_{0}-C_{T+1} + 2T\left\vert \mathcal{A}-1\right\vert \mathcal{O}(\beta ^{2})  \\
& = \sum\limits_{a\in \mathcal{A}}[C^*-(\tilde{C}_{-a}+C_{a,0})+%
\min_{a}\{C_{a,0}+\tilde{C}_{-a}-C^{\ast }\}] \nonumber\\&\hspace{6pt}+\sum\limits_{a\in \mathcal{A}%
}(r_{a,T+1}-r_{a,0})+C_{0}-C^*+ \delta_T(\beta)
\label{eq_budget_T_prime}
\end{align}
where 
\begin{align}
\delta _{T}(\beta) &:=\sum\limits_{a\in \mathcal{A}}[C_{T+1}-C^{\ast }-(\tilde{%
C}_{-a,T+1}-\tilde{C}_{-a})]\nonumber\\& \hspace{6pt}+C^{\ast }-C_{T+1}+ 2T\left\vert \mathcal{A}-1\right\vert \mathcal{O}(\beta ^{2}).
\end{align}We get \eqref{eq_budget_T_prime} by adding and subtracting $\sum_{a\in\ccalA} C^*$, $\sum_{a\in\ccalA} \tilde C_{-a}$, and $C^*$. Now we use \eqref{efficient_coupling_prime} to get
\begin{align}
   \widehat{B}(T,\beta)&\geq \sum_{a\in\ccalA} [r_{a,T+1}- r_{a,0}]+\delta _{T}(\beta).
\end{align}Setting initial locational marginal prices and capacity prices associated with the tielines to zero, we have $r_{a,0}=0$. We have $\delta_T(\beta) \rightarrow 0$ as $T\to \infty$ and thus we get 
\begin{align}
\lim_{T\rightarrow \infty,\beta \rightarrow 0^+ } \widehat{B}(T,\beta) &\geq \sum\limits_{a\in \mathcal{A}%
}\big[\sum_{(i,j)\in \mathcal{T}_{a}}{\alpha _{a^{\prime },j}^{\ast }}{T}%
_{a,(i,j)}^{\ast }-\sum_{(i,j)\in \mathcal{T}_{a}}\frac{\mu _{(i,j)}^{\ast }%
}{2}(\left\vert {T}_{a,(i,j)}^{\ast }\right\vert -\bar{T}_{a,(i,j)})\big] \nonumber
\end{align}%
The result follows by noting that ${T}_{a,(i,j)}^{\ast }=-{T}%
_{a^{\prime },(j,i)}^{\ast }$ and the fact the the capacity price associated terms above converge to zero by complementary slackness. 

\section{Acknowledgements}
Alfredo Garcia gratefully acknowledges partial support from NSF grant $\#1829552$ and the Texas A\&M Energy Institute. Ceyhun Eksin and Furkan Sezer were, in part, supported by grants from the National Science Foundation (NSF CCF-2008855 and NSF ECCS-1953694).

\bibliographystyle{IEEEtran}
\bibliography{ref}
\vspace{-5pt}

\end{document}